\documentclass[runningheads, numbook, envcountsame]{amsart}
\usepackage{amssymb, amsmath, amscd, amsthm}
\usepackage{xypic, rotating}
\usepackage{graphicx}
\usepackage{epsf, mathptmx}
\input xy
\xyoption{all}

\newcommand{\R}{\ensuremath{\mathbb{R}}}
\newcommand{\F}{\ensuremath{\mathbb{F}}}
\newcommand{\Z}{\ensuremath{\mathbb{Z}}}
\newcommand{\T}{\ensuremath{\mathbb{T}}}

\newcommand{\sC}{\ensuremath{{\mathcal C}}}
\newcommand{\sD}{\ensuremath{{\mathcal D}}}
\newcommand{\sE}{\ensuremath{{\mathcal E}}}

\newcommand{\sM}{\ensuremath{{\mathcal M}}}
\newcommand{\W}{\ensuremath{{\mathcal W}}}
\newcommand{\Mbar}{\ensuremath{\overline{\mathcal M}}}
\newcommand{\sO}{\ensuremath{{\mathcal O}}}

\newcommand{\C}{\ensuremath{\mathbb{C}}}

\newcommand{\g}{\ensuremath{\mathfrak{g}}}

\newcommand{\Hom}{{\rm Hom}}

\newcommand{\Tor}{{\rm Tor}}
\newcommand{\im}{{\rm im \,}}
\newcommand{\ungamma}{\underline{\gamma}}
\newcommand{\ovgamma}{\overline{\gamma}}
\newcommand{\unrho}{\underline{\rho}}
\newcommand{\uncirc}{\underline{\circ}}
\newcommand{\ovrho}{\overline{\rho}}
\newcommand{\unG}{\underline{G}}
\newcommand{\sma}{\land}

\begin{document}
\title{Equivariant operads, string topology, and Tate cohomology}
\author{Craig Westerland}
\address{University of Wisconsin-Madison \\
Mathematics Department \\
480 Lincoln Dr. \\
Madison, WI 53706-1388, USA}


\newtheorem{theorem}{Theorem}[section]
\newtheorem{proposition}[theorem]{Proposition}
\newtheorem{convention}[theorem]{Convention}
\newtheorem{comment}[theorem]{Comment}
\newtheorem{lemma}[theorem]{Lemma}
\newtheorem{obs}[theorem]{Observation}
\newtheorem{conjecture}[theorem]{Conjecture}
\newtheorem{corollary}[theorem]{Corollary}
\newtheorem{claim}[theorem]{Claim}
\newtheorem{construction}[theorem]{Construction}
\theoremstyle{definition}
\newtheorem{definition}[theorem]{Definition}
\theoremstyle{remark}
\newtheorem{remark}[theorem]{Remark}
\newtheorem{example}[theorem]{Example}
%

\bibliographystyle{amsalpha}
\begin{abstract}
From an operad $\sC$ with an action of a group $G$, we construct new operads using the homotopy fixed point and orbit spectra.  These new operads are shown to be equivalent when the generalized $G$-Tate cohomology of $\sC$ is trivial.  Applying this theory to the little disk operad $\sC_2$ (which is an $S^1$-operad) we obtain variations on Getzler's gravity operad, which we show governs the Chas-Sullivan string bracket. 

\end{abstract}
\maketitle
%
%
%
%
%
%
\section{Introduction} \label{intro}

In their foundational paper \cite{cs}, Chas and Sullivan constructed a variety of algebraic structures on the singular and $S^1$-equivariant homology of the free loop space $LM = Map(S^1, M)$ of a closed $d$-manifold $M$.  In the beautiful \cite{cj}, Cohen and Jones showed that the Chas-Sullivan operations on $H_*(LM)$ are governed by an operad -- the cactus operad, or equivalently, the two-dimensional framed little disk operad.  One purpose of this paper is to do likewise for the Chas-Sullivan operations in the Borel equivariant homology, $H_*^{S^1}(LM) = H_*(ES^1 \times_{S^1} LM)$.

Recall that Chas and Sullivan introduced the \emph{string bracket}, a graded Lie bracket of dimension $2-d$:
$$[ \cdot, \cdot ] : H_p^{S^1} (LM) \otimes H_q^{S^1} (LM) \to H_{p+q +2-d}^{S^1} (LM)$$
This was defined as follows: for classes $a, b \in H_*^{S^1}(LM)$, 
$$[a, b] := p(\tau(a) \cdot \tau(b))$$
where $\cdot$ is the Chas-Sullivan \emph{loop product} on $H_*(LM)$, $\tau$ is the $S^1$-transfer 
$$\tau:H_*^{S^1}(LM) \to H_{*+1}(LM)$$
and $p$ is the projection to the quotient, $H_*(LM) \to H_*^{S^1}(LM)$.

It is remarkable that bracketing a ring multiplication with $\tau$ and $p$ produces a Lie bracket.  One is led to wonder whether this is an example of a construction of a more general nature.  Furthermore, Chas and Sullivan define a family of $k$-ary operations $\overline{m}_k: H_*^{S^1}(LM)^{\otimes k} \to H_*^{S^1}(LM)$ by
$$\overline{m}_k(a_1 \otimes \dots \otimes a_k) = p(\tau(a_1) \cdots \tau(a_k))$$
From an operadic point of view, this construction of operations of higher ``arity'' is very natural.  The main application of this paper to string homology will be to give the action of an operad which governs these operations.  

In \cite{getz2d,getz0}, Getzler defined the \emph{gravity operad} $Grav$ in the category of graded groups using the (open) moduli spaces of points in $\C P^1$:
$$Grav(k):= \Sigma H_*(\sM_{0, k+1}).$$
In the language of \cite{getz0}, this is actually the operadic desuspension $\Lambda^{-1} Grav$.  For clarity we will forego the desuspension in this paper.  Unlike many familiar operads (such as the associative, commutative, and Gerstenhaber operads), the gravity operad is not generated by a finite number of operations.  However, in \cite{getz2d}, Getzler defines an infinite family of degree $1$ operations corresponding to the generator of $H_0(\sM_{0, k+1}) = \Z$:
$$\{ a_1, \dots, a_k \} \in Grav(k), \; k \geq 2$$
This is a single $k$-ary operation, denoted by braces; the $a_i$ are dummy variables.  Collectively these brackets generate $Grav$, subject to a generalized Jacobi relation:
\begin{multline*}
\{ \{ a_1, \dots, a_k \}, b_1, \dots, b_l \} = \\
\sum_{1 \leq i < j \leq k} (-1)^{\epsilon(i, j)} \{ \{ a_i, a_j \}, a_1, \dots, \hat{a}_i, \dots, \hat{a}_j, \dots, a_k, b_1, \dots b_l \}
\end{multline*}
where $\epsilon(i, j) = (|a_1| + \dots + |a_{i-1}|)|a_i| + (|a_1| + \dots + |a_{j-1}|)|a_j| + |a_i||a_j|$.  This relation makes the binary bracket into a Lie bracket (of graded degree $1$).

\begin{theorem} \label{string_homology_theorem}

Let $h_*$ be a ring homology theory which supports an orientation for $TM$.  Then the (generalized, shifted) string homology $\Sigma^{1-d} h_*(LM_{hS^1})$ of a $d$-manifold is an algebra over the generalized gravity operad 
$$Grav_{h_*} := \{ \Sigma h_*(\sM_{0, k+1}) \}$$
Moreover, when $h_*$ is taken to be singular homology, the brackets encode the Chas-Sullivan operations:
$$\overline{m}_k(a_1 \otimes \dots \otimes a_k) = \{ a_1, \dots, a_k \}$$
for $a_i \in H_*(LM_{hS^1})$.

\end{theorem}

\noindent Here, for a graded group $A_*$ and an integer $s$, we use the notation $\Sigma^s A_*$ to mean the graded group whose $n+s^{\rm th}$ term is $A_n$.

We will give a geometric construction of the gravity operad from an equivariant stable homotopy theory point of view, much along the lines of Getzler's construction in \cite{getz2d}.  In fact, this ``transfer operad'' will be constructed as an equivariant version of the cactus operad used in \cite{cj}, and its action will be induced by Cohen and Jones's construction.

More generally, the purpose of this paper is to study how operads and their algebras, when equipped with a group action, give rise to new operads and algebras after applying various constructions in the equivariant category.  From a $G$-equivariant operad $\sC$ we will define a new (non-equivariant) operad in the stable category, whose terms are constructed as certain Thom spaces over $\sC(k)_{hG}$, and whose substitution map employs the $G$-transfer map.  For the little disk operad $\sC_2$ (with an action of $S^1$), this produces an incarnation of the gravity operad.   Equivariant algebras $X$ over the original operad give rise to algebras over the new operad by the same construction -- a Thom space over $X_{hG}$.  Very roughly, this gives Theorem \ref{string_homology_theorem} when $X = LM$ using Cohen and Jones's theorem.  There is similarly an analogue for homotopy fixed point spectra, and one may compare the two constructions using Tate cohomology.

\section{Summary of the paper}

We begin with general facts about constructing new operads and algebras from equivariant ones, and adapt these constructions to the string homology context.  Our constructions take us into the equivariant stable category, but one can ask similar questions in the category of chain complexes.  As should be obvious to even a casual reader of \cite{getz2d,getz0}, when the group in question is $S^1$, chain complex analogues of these constructions were known to Getzler over a decade ago; this is most apparent in Corollary \ref{homology_corollary}.

\subsection{Constructing new operads}

Though they are most commonly studied in the categories of topological spaces and (differential) graded vector spaces, operads can be defined in any symmetric monoidal category.  For instance, for a group $G$, an operad $\sC$ in the symmetric monoidal category of $G$-spaces is called a \emph{$G$-operad}; this is distinguished from usual operads in spaces in that each $\sC(k)$ is a $G$-space, and the action of $\Sigma_k$ on $\sC(k)$ and the substitution map 
$$\gamma: \sC(k) \times \sC(n_1) \times \dots \times \sC(n_k) \to \sC(\sum_{i=1}^k n_i)$$
are $G$-equivariant (here we employ the diagonal $G$ action on the domain of $\gamma$).  For a discussion of such gadgets, we refer the reader to \cite{salvwahl}.  An excellent example is the $k$-dimensional little disk operad $\sC_k$, which admits an action of $SO(k)$ by rigid rotations within $\R^k$.

\begin{figure}
\centering
\includegraphics{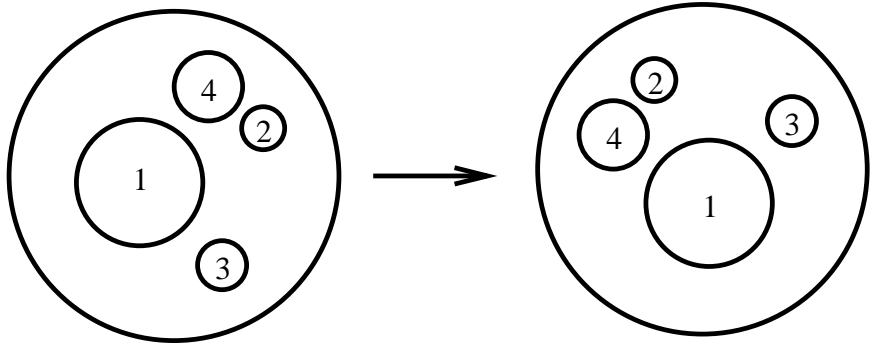}
\caption{The action of $i \in S^1 = SO(2)$ on an element of $\sC_2(4)$.}
\label{disc_fig}
\end{figure} 

In this paper we will also consider operads in certain categories of spectra.  Here, by the category of spectra, we mean a modern point-set category of spectra, specifically the $S$-modules of \cite{EKMM}.  We will also employ the category of (genuine) $G$-spectra, indexed on a fixed complete $G$-universe $U$, as well as the category of naive $G$-spectra, indexed on the $G$-trivial universe $U^G$.  For an introduction to these categories, we refer the reader to \cite{lms,manmay}.

An operad in the category of spectra is a collection of $S$-modules $\sC(k)$ ($k>0$) equipped with an action of $\Sigma_k$ on $\sC(k)$ and a substitution map
$$\gamma: \sC(k) \sma \sC(n_1) \sma \dots \sma \sC(n_k) \to \sC(\sum_{i=1}^k n_i)$$
satisfying the appropriate axioms (see, e.g., \cite{may}).  We will \emph{not} insist that our operads be unital.  An operad $\sC$ in $G$-spectra or naive $G$-spectra is given by the same data, where $\sC(k)$ are either $G$-spectra or naive $G$-spectra, and $\gamma$ and the $\Sigma_k$ action are appropriately equivariant.  We will refer to such operads as \emph{stable $G$-operads} and \emph{naive $G$-operads}, respectively.  The following simple fact is one of the best methods of producing examples of naive $G$-operads.

\begin{proposition}

If $\sC$ is a $G$-operad in the category of spaces, then the suspension spectra $\Sigma^{\infty} \sC(k)_+$ form a naive $G$-operad, $\Sigma^{\infty} \sC_+$.  Furthermore, if $X$ is an algebra over $\sC$ in the category of $G$-spaces, then $\Sigma^{\infty} X_+$ is an algebra over $\Sigma^{\infty} \sC_+$.

\end{proposition}

We find the following diagram of functors amongst the three categories at hand enlightening; it summarizes most of the constructions that we will make.  Note however, that it does \emph{not} commute; in some sense Tate cohomology measures its failure to commute.

\begingroup
\small
$$\xymatrixcolsep{3pc}\xymatrix{
 & naive \; G-spectra \ar[r]^-{i_*} & G-spectra \ar[d]^-{i^*} & \\
G-spaces \ar[r]_-{\Sigma^{\infty} \cdot_+} & naive \; G-spectra \ar[r]_-{F(EG_+, \cdot )} \ar[u]^-{EG_+ \sma \cdot} & naive \; G-spectra \ar[r]_-{\cdot^G} & spectra
}$$
\endgroup
Here $i^*$ and $i_*$ are change of universe functors; they are induced by the inclusion $i:U^G \to U$.  The forgetful functor $i^*$ produces a spectrum with the same spaces, but indexed on the sub-universe $U^G$, and $i_*$ builds in non-trivial representations.  The space $EG$ is a free, contractible $G$-space, and $EG_+ \sma \cdot$ and $F(EG_+, \cdot)$ are smash product and function spectra.  Finally, $\cdot^G$ denotes the fixed point spectrum functor.

Let us examine some of the possible composites of the functors in this diagram.  For a naive $G$-spectrum $E$, $E^{hG} := F(EG_+, E)^G$ is called the \emph{homotopy fixed point spectrum}.  The spectrum $E_{hG} := (EG_+ \sma E) / G$ is called the \emph{homotopy orbit spectrum} or \emph{Borel construction}.  If $G$ is assumed to be compact Lie (see \cite{klein_dual,rognes} for ways to relax this assumption), the transfer map defines an equivalence
$$\tau^G: (EG_+ \sma E \sma S^{Ad_G})/G \to (i^* i_*(EG_+ \sma E))^G$$
where $S^{Ad_G}$ is the one-point compactification of the Lie algebra of $G$, equipped with the adjoint action of $G$ (see for instance \cite{gm} for a proof).  For brevity, we will write $E_{bG} := (i^* i_*(EG_+ \sma E))^G$.  We see that $E_{bG}$ is closely related to the Borel construction of the $G$-action on $E$.

A short summary of the first half of this paper is that all of the functors in the diagram above carry operads in the domain category to operads in the target category.  From our point of view the most fundamental is the fixed point construction.  Fix a naive $G$-operad $\sC$, as well as a $\sC$-algebra $X$ in the category of naive $G$-spectra.

\begin{theorem} \label{fixed_op_theorem}

The fixed point spectra $\sC(k)^G$ form a operad $\sC^G$ in spectra.  Furthermore, $X^G$ is a $\sC^G$-algebra.  We call $\sC^G$ the \emph{fixed point operad}.

\end{theorem}

\begin{remark} 

A way to reformulate the condition that $X$ is a $\sC$-algebra in naive $G$-spectra is as follows. Recall from \cite{salvwahl} the notion of a semidirect product 
$$(\sC \rtimes G)(k) = \sC(k) \times G^{\times k}$$
The substitution map of $\sC \rtimes G$ is that of $\sC$, twisted by the action of $G$:  
\begin{multline*}
(c; g_1, \dots, g_k) \circ ((d_1; h_1^1, \dots ,h_{n_1}^1), \dots , (d_k; h_1^k, \dots , h_{n_k}^k)) = \\
(c \circ (g_1 d_1, \dots, g_k d_k); g_1 h_1^1, \dots, g_k h_{n_k}^k)
\end{multline*}
For instance, the semidirect product $\sC_k \rtimes SO(k)$ is the $k$-dimensional \emph{framed} little disk operad.

Salvatore and Wahl show that a $\sC \rtimes G$-algebra is simply a $\sC$-algebra in the category of $G$-spaces.  For a naive $G$-operad $\sC$, we may similarly form the semidirect product operad $\sC \rtimes G$, and a $\sC \rtimes G$-algebra is the same thing as a $\sC$-algebra in the category of naive $G$-spectra.

\end{remark}

The following is the main application of Theorem \ref{fixed_op_theorem}.

\begin{corollary} \label{homotopy_corollary}

The spectra $\sC(k)^{hG}$ (respectively $\sC(k)_{bG}$) form an operad $\sC^{hG}$ (resp. $\sC_{bG}$) in spectra, and $X^{hG}$ (resp. $X_{bG}$) is an algebra over $\sC^{hG}$ (resp. $\sC_{bG}$).

\end{corollary}

\noindent We will call $\sC^{hG}$ the \emph{homotopy fixed point operad}, and  $\sC_{bG}$ the \emph{transfer operad}.  If $\sC$ is unital and the unit is fixed by the $G$ action, then $\sC^{hG}$ is unital as well, but $\sC_{bG}$ is generally not.

Notice that if $E = \Sigma^{\infty} X_+$ is the suspension spectrum of a space, then
$$E_{bG} \simeq  (EG_+ \sma E \sma S^{Ad_G})/G = \Sigma^{\infty} (X_+ \sma S^{Ad_G})_{hG}$$
is the suspension spectrum of the Thom space $X_{hG}^{Ad_G}$ of the flat $G$-bundle over $X_{hG}$ whose fiber is the adjoint representation of $G$.  So for any ring homology theory $h_*$ for which the vector bundle $Ad_G$ admits a Thom isomorphism, 
$$h_*(X_{bG}) \cong \Sigma^{\dim (G)} h_*(X_{hG})$$
Thus Corollary \ref{homotopy_corollary} gives the following:

\begin{corollary} \label{homology_corollary}

For a $G$-operad $\sC$ in spaces and a $\sC \rtimes G$-algebra $X$, the homology of the transfer operad
$$\Sigma^{\dim(G)} h_*(\sC_{hG})$$
is an operad in graded $h_*$-modules, and $\Sigma^{\dim(G)} h_*(X_{hG})$ is an algebra over it.

\end{corollary}

\subsection{Examples of new operads}

The first example of an equivariant operad to study is the 2-dimensional little disks operad $\sC_2$, which is an $S^1 = SO(2)$-operad.

\begin{proposition} \label{gravity_prop}

The homology of the transfer operad $(H_*({\sC_2}_{bS^1}))_{>1}$ is isomorphic to the gravity operad $Grav$ (away from the unary term).

\end{proposition}

Essentially, this is a rephrasing of some of the constructions of \cite{getz2d} in the category of spectra.  Interestingly, the unary part of this operad can be thought of in moduli-theoretic terms as well; see section \ref{2d_section}.  We will prove this proposition in that section, but we point out here that since $\sC_2(2)$ is equivariantly homotopy equivalent to $S^1$, 
$$
H_k(\sC_2(2)_{bS^1}) = H_k(\Sigma \Sigma^{\infty} ({\sC_2(2)}_{hS^1})_+) = H_{k-1}(pt.) = \left\{ \begin{array}{ccc}
                                                                 \Z & , & k = 1 \\
                                                                 0  & , & k \neq 1
                                                                 \end{array}
                                                           \right.
$$
This sole generator defines the binary Lie bracket of the gravity operad, and gives rise to the Chas-Sullivan string bracket in the string homology context.  In fact, it induces an embedding of the (graded) Lie operad as a suboperad of the homology operad $H_*({\sC_2}_{bS^1})$.  

We study a four-dimensional analogue of $Grav$ in section \ref{4d_section} using the transfer operad fror the $SU(2)$-action on $\sC_4$.

For a group $G$ one may form an operad $\unG$ with $\unG(k) = G^{\times k}$.  $\unG$ is a $G$-operad via (termwise) conjugation.  The resulting transfer and homotopy fixed point operads are related to the string topology of the classifying space $BG$; they are studied in section \ref{loop_subsection}.

\subsection{Tate cohomology}


Recall that, for any naive $G$-spectrum $E$, there is a \emph{norm map}
$$n^{G}: E_{bG} \to E^{hG}$$
whose cofiber is the generalized Tate cohomology spectrum $E^{tG}$.  We use this definition of Tate cohomology to relate the homotopy fixed point and transfer operads:

\begin{proposition} \label{tate_proposition}

The norm map
$$n^{G}: \sC_{bG} \to \sC^{hG}$$
is a map of operads, and if the generalized Tate cohomology spectra $\sC(k)^{tG}$ are trivial, it is a weak equivalence of operads.

\end{proposition}

A simple application of this proposition is the following:

\begin{corollary} \label{c2_corollary}

There is a weak equivalence of operads $({\sC_2}_{bS^1})_{>1} \simeq(\sC_2^{hS^1})_{>1}$.

\end{corollary}

Here $(\sE)_{>1}$ denotes the operad whose first term is trivial, and whose higher terms are those of $\sE$: $(\sE)_{>1}(k)  = \sE(k)$ if $k>1$.  As a result, we see that 
$$(H_*(\sC_2^{hS^1}))_{>1} \cong Grav$$

\subsection{String homology}

Recall Cohen and Jones' result \cite{cj} that the homology of the cactus operad governs the string topology operations on $H_*(LM)$.  It is a theorem of Voronov's \cite{voruniv} that the cactus operad is equivalent to the framed disk operad $\sC_2 \rtimes S^1$.  So in light of Proposition \ref{gravity_prop} and Corollary \ref{homology_corollary} (applied to the case $G = S^1$, $\sC = \sC_2$, and $X = LM$), Theorem \ref{string_homology_theorem} should not be surprising.

Some substantial modifications are necessary, however.  Recall, for instance, that Cohen and Jones do not prove in \cite{cj} that $LM$ is an algebra over the cactus operad.  Rather, they give maps that induce an algebra structure after application of any generalized homology theory $h_*$ that supports an orientation of $TM$.  We get around this difficulty in the equivariant case in the same way they did, via certain Pontrjagin-Thom collapse maps.  Additionally, we have to modify the construction of the transfer operad for the cactus operad, whose equivariance is not as straightforward as the little disk operad.  

We also obtain applications to a homotopy fixed point variation of string homology:

\begin{corollary} \label{fixed_string_corollary}

Let $h_*$ be a ring homology theory which supports an orientation for $TM$.  Then $\Sigma^{-d} h_*(LM^{hS^1})$ is an algebra over the $h_*$-homology of the homotopy fixed point operad, $h_*({\sC_2}^{hS^1})$.

\end{corollary}

Consequently $\Sigma^{-d} H_*(LM^{hS^1})$ is similarly a gravity algebra, and in particular admits a Lie bracket analogous to the Chas-Sullivan string bracket.  We study this in detail in section \ref{string_homology_section}, and prove a splitting theorem for this structure on a ``continuous'' version of this homology.

%
%

\subsection{Structure of the paper}

In section \ref{new_op_section} we will construct the fixed point, homotopy fixed point, and transfer operads.  Additionally we discuss a simplification, fixed point operads in the category of spaces.  Section \ref{string_homology_section} is devoted to adapting these results to the string topology context.  Beforehand, however, in section \ref{tate_section} we explore the relationship with Tate cohomology and prove Proposition \ref{tate_proposition}.  In section \ref{dual_section} we study the (continuous) homology of homotopy fixed point operads.  Section \ref{example_section} is devoted to a detailed study of the homology of the transfer and homotopy fixed point operads for the $S^1$-action on $\sC_2$, the $SU(2)$-action on $\sC_4$, and the $G$-action on $\unG$.

\section{Constructing new operads} \label{new_op_section} 

In this section we prove Theorem \ref{fixed_op_theorem} and deduce Corollary \ref{homotopy_corollary}.  This will require some careful point-set manipulations of spectra; a more concrete (but less precise) construction of the transfer operad is given briefly in section \ref{sketch_subsection}.  The construction there is sufficient to show that the transfer operad is an operad in the stable homotopy category.  The reader who is only interested in this level of detail is encouraged to skip sections \ref{fixed_op_section} and \ref{homotopy_section} where these results are promoted to the category of $S$-modules.

\subsection{Operads} \label{defn_subsection}

Operads were introduced by May in \cite{may}.  There are now at least three ways to define them.  We shall employ all of these definitions in the sequel.  We briefly review these definitions in this section; a detailed exposition is given in \cite{may_defn}.

In May's original definition, an operad $\sC$ in a symmetric monoidal category $(C, \otimes)$ with finite coproducts is defined as a collection of objects $\sC(k)$ for $k \geq 1$, equipped with an action of $\Sigma_k$ on $\sC(k)$ and a collection of substitution maps
$$\gamma: \sC(k) \otimes \sC(n_1) \otimes \dots \otimes \sC(n_k) \to \sC(\sum n_i)$$
which are appropriately associative and equivariant.

In \cite{kelly,gj,ching}, operads in $C$ were presented as monoids in the category of symmetric sequences in $C$, which is made into a monoidal category using the composition product of symmetric sequences.  This description neatly packages the axioms of associativity and equivariance in the previous definition into the construction of this monoid.

Let $P$ be the category whose objects are the positive integers, with $\Hom_P(m, n)$ empty if $m \neq n$, and the symmetric group $\Sigma_m$ if $m=n$.  Recall that a \emph{symmetric sequence in $C$} is a functor $\sC: P \to C$.  The \emph{composition product} $\sC \circ \sD$ of two such sequences $\sC$, $\sD$ is defined by
$$\sC \circ \sD(m) := \bigsqcup Ind (\sC(k) \otimes \sD(n_1) \otimes \dots \otimes \sD(n_k))$$
The coproduct is taken over all $k$ and all partitions $n_1 + \dots + n_k = m$, and $Ind$ is induction on the inclusion $\Sigma_k \times \Sigma_{n_1} \times \dots \times \Sigma_{n_k} \leq \Sigma_m$:
$$Ind (\sC(k) \otimes \sD(n_1) \otimes \dots \sD(n_k)) := \frac{\Sigma_m \times (\sC(k) \otimes \sD(n_1) \otimes \dots \otimes \sD(n_k))}{\Sigma_k \times \Sigma_{n_1} \times \dots \times \Sigma_{n_k}} $$
Here, for a set $S$ and an object $A$ of $C$, $S \times A$ is the coproduct $\sqcup_S A$.

The composition product makes the category of symmetric sequences (that is, the functor category $[P, C]$) into a (nonsymmetric) monoidal category.  An operad is then simply a (not necessarily unital) monoid for $\circ$ in this category; it is equipped with an associative map 
$$\mu: \sC \circ \sC \to \sC$$
Comparing this to the first definition, $\mu = \sqcup Ind(\gamma)$.

Another definition employs the partial substitutions, or $\circ_i$ products.  In this formulation, an operad is a symmetric sequence $\sC$ equipped with maps
$$\circ_i: \sC(m) \otimes \sC(n) \to \sC(m+n-1), \; 1 \leq i \leq m$$
which are associative and equivariant with respect to the group $\Sigma_{m-1} \times \Sigma_n$.  Here $\Sigma_{m-1} \times \Sigma_n$ acts on the domain of $\circ_i$ as the subgroup of $\Sigma_m \times \Sigma_n$ which stabilizes $i \in \{ 1, \dots, m \}$, and on the range as the subgroup of $\Sigma_{m+n-1}$ that does not permute the first $m-1$ elements of $\{1, \dots, m+n-1 \}$ with the last $n$ elements.

If we require our operads to be unital, all of these definitions agree: 
$$c \circ_i d = \gamma(c; 1, \dots, 1, d, 1, \dots, 1),$$
where $d$ is placed in the $i^{\rm th}$ position .  However, if we do not assume unitality, they are inequivalent: the $\circ_i$ definition is stronger than the classical definition.  In this paper, we will employ the weaker definition.

\subsection{A naive construction of the transfer operad} \label{sketch_subsection}

If we begin with a compact Lie group $G$, and a $G$-operad $\sC$ in topological spaces, there is an instructive construction of the substitution map for the transfer operad $(\Sigma^{\infty} \sC_+)_{bG}$ using relative transfer maps.  Since transfer maps are only well-defined up to homotopy (requiring many choices), this produces an operad up to homotopy.  

For a $G$-space $X$, one may define a vector bundle $Ad_G$ over $X_{hG}$ as
$$(X \times EG) \times_G \g \to (X \times EG) / G = X_{hG}.$$
Here the fibre is $\g$, the Lie algebra of $G$, equipped with the adjoint action of $G$.  We will denote the Thom space of this bundle $X_{hG}^{Ad_G}$.  Recall from \cite{madsch} that for closed subgroups $H<G$ there is a relative transfer map:
$$\tau_H^G: \Sigma^{\infty} X_{hG}^{Ad_G} \to \Sigma^{\infty} X_{hH}^{Ad_H}$$

We will define the structure of an operad on the collection 
$$\sC_{hG}^{Ad_G} = \{ \Sigma^{\infty} \sC(k)_{hG}^{Ad_G}, \, k \geq 1 \}$$
as follows: because $\sC$ is a $G$-operad, the substitution in $\sC$
$$\gamma: \sC(k) \times \sC(n_1) \times \dots \times \sC(n_k) \to \sC(\sum_i n_i)$$
is $G$-equivariant (where $G$ acts diagonally on the left).  Therefore $\gamma$ descends to a map
$$\gamma_{hG}: (\sC(k) \times \sC(n_1) \times \dots \times \sC(n_k) )_{hG} \to \sC(\sum_i n_i)_{hG} \leqno{(*)} $$
This is covered by a map of bundles $Ad_G \to Ad_G$ and so extends to the Thom spaces of $Ad_G$.

The space $\sC(k) \times \sC(n_1) \times \dots \times \sC(n_k)$ is acted upon by $G^{\times k+1}$, and we can think of the left hand side of $(*)$ as the homotopy quotient by the diagonal subgroup $\Delta(G) < G^{\times k+1}$.  The relative transfer map $\tau_{\Delta(G)}^{G^{\times k+1}}$ in this case is
$$\Sigma^{\infty} (\sC(k) \times \sC(n_1) \times \dots \times \sC(n_k))_{hG^{\times k+1}}^{Ad_{G^{\times k+1}}} \to \Sigma^{\infty} (\sC(k) \times \sC(n_1) \times \dots \times \sC(n_k) )_{hG}^{Ad_G}$$

The left hand side is equivalent to
$$\Sigma^{\infty} \sC(k)_{hG}^{Ad_G} \sma \Sigma^{\infty} \sC(n_1)_{hG}^{Ad_G} \sma \dots \sma \Sigma^{\infty} \sC(n_k)_{hG}^{Ad_G}$$
so the composition $\ungamma := (\Sigma^{\infty} \gamma_{hG}^{Ad_G}) \circ \tau_{\Delta(G)}^{G^{\times k+1}}$ is a map
$$\ungamma: \Sigma^{\infty} \sC(k)_{hG}^{Ad_G} \sma \Sigma^{\infty} \sC(n_1)_{hG}^{Ad_G} \sma \dots \sma \Sigma^{\infty} \sC(n_k)_{hG}^{Ad_G} \to \Sigma^{\infty} \sC(\sum_i n_i)_{hG}^{Ad_G}$$
This will serve as a substitution map for the operad $\sC_{hG}^{Ad_G}$.  

One may similarly obtain a partial substitution 
$$\uncirc_i: \Sigma^{\infty} \sC(k)_{hG}^{Ad_G} \sma \Sigma^{\infty} \sC(l)_{hG}^{Ad_G}  \to \Sigma^{\infty} \sC(k+l-1)_{hG}^{Ad_G}$$
as the composite $\uncirc_i = (\circ_i)_{hG} \circ \tau_{\Delta(G)}^{G \times G}$.

In the next sections we will construct a \emph{strict} operad substitution for the collection of spectra $\sD_{bG} = \{ \sD(k)_{bG}, \, k\geq 1 \}$ (for naive $G$-operads $\sD$).  This is related to the construction above: after conjugating by the equivalence
$$(\Sigma^\infty \sC(k)_+)_{bG} \simeq \Sigma^{\infty} \sC(k)^{Ad_G}_{hG},$$
$\ungamma$ is homotopic to the substitution map for $(\Sigma^{\infty} \sC_+)_{bG}$.

\subsection{The fixed point operad} \label{fixed_op_section}

The heart of the proof of Theorem \ref{fixed_op_theorem} is the construction of the operad substitution map for $\sC^G$ from $\sC$.  This relies upon the following fact:

\begin{lemma} \label{monoidal_lemma}

Let $X$ and $Y$ be two naive $G$-spectra.  Then there is a natural transformation 
$$i_{X, Y}: X^G \sma Y^G \to (X \sma Y)^G$$
where, on the right, $G$ acts diagonally on $X \sma Y$.  This is associative in the sense that for any $X$, $Y$, and $Z$, the following diagram commutes
$$
\xymatrix{
X^G \sma Y^G \sma Z^G \ar[r]^-{i_{X, Y} \sma 1} \ar[d]_-{1 \sma i_{Y, Z}} & (X \sma Y)^G \sma Z^G \ar[d]^-{i_{X \sma Y, Z}} \\
X^G \sma (Y \sma Z)^G \ar[r]_-{i_{X, Y \sma Z}}  & (X \sma Y \sma Z)^G
}.
$$
Furthermore, $i$ is symmetric, in that the following diagram also commutes:
$$\xymatrix{
X^G \sma Y^G \ar[r]^-{i_{X, Y}} \ar[d]_T & (X \sma Y)^G \ar[d]^{T^G} \\
Y^G \sma X^G \ar[r]^-{i_{Y, X}} & (Y \sma X)^G
}$$
where $T$ is the symmetry isomorphism.

\end{lemma}

This natural transformation is studied in section VI.3 of \cite{manmay}.  Such a functor is called a ``lax'' symmetric monoidal functor.  Thanks to the first commutative diagram, we can speak unambiguously of \emph{the} map
$$i: X_1^G \sma \dots \sma X_n^G \to (X_1 \sma \dots \sma X_n)^G$$
for any family $X_1, \dots, X_n$ of naive $G$-spectra.

\begin{definition}

For a naive $G$-operad $\sC$ with substitution map $\gamma$, the operad substitution for $\sC^G$ is defined as the composite $\ovgamma :=\gamma^G \circ i$:
$$\sC(k)^G \sma \sC(n_1)^G \sma \dots \sma \sC(n_k)^G \to (\sC(k) \sma \sC(n_1) \sma \dots \sma \sC(n_k))^G \to \sC(\sum n_i)^G$$

\end{definition}

To see that this actually gives $\sC^G$ the structure of an operad, we need the following proposition, which is a mild generalization of a result of \cite{kelly}.

\begin{proposition} \label{lax_prop}

Let $C$ and $D$ be symmetric monoidal categories, and $\sC$ an operad in $C$.  If $F: C \to D$ is a lax symmetric monoidal functor, then the collection $F\sC = \{ F\sC(k), \, k \geq 1 \}$ forms an operad in $D$.

\end{proposition}

Theorem \ref{fixed_op_theorem} follows as a corollary of this proposition and Lemma \ref{monoidal_lemma}.


\begin{proof}

Postcomposition with $F$ defines a functor
$$F_* : [P, C] \to [P, D]$$
In other words, $F$ carries a symmetric sequence in $C$ to a symmetric sequence in $D$.  Further, if $F$ is lax symmetric monoidal, then $F_*$ is lax monoidal (with respect to the composition product of symmetric sequences).  Thus $F_*$ carries an operad in $C$ (a monoid in $[P, C]$) to an operad in $D$, since lax monoidal functors carry monoids to monoids.  Concretely, the operad substitution in $F\sC$ is given by the composite $(F \gamma) \circ i$:
$$F\sC(k) \sma F\sC(n_1) \sma \dots \sma F\sC(n_k) \to F(\sC(k) \sma \sC(n_1) \sma \dots \sma \sC(n_k)) \to F\sC(\sum n_i)$$

\end{proof}

\subsection{The transfer and homotopy fixed point operads} \label{homotopy_section}

Using the results of the previous section, we will prove Corollary \ref{homotopy_corollary}.  We begin with the following fairly trivial constructions.

\begin{lemma} \label{X_lemma}

If $\sC$ is an operad in spectra and $X$ is a space, then the collection $X_+ \sma \sC$ of spectra
$$\{ X_+ \sma \sC(k) \}_{k \geq 1}$$
forms an operad in spectra.  Likewise, the collection $F(X_+, \; \sC)$
$$\{ F(X_+, \; \sC(k)) \}_{k \geq 1}$$
forms an operad.  Further, if $X$ is a $G$-space and $\sC$ a naive $G$-operad, then both $X_+ \sma \sC$ and $F(X_+, \; \sC)$ are naive $G$-operads.

\end{lemma}

The $G$-action on $X_+ \sma \sC(k)$ is diagonal, and on $F(X_+, \; \sC(k))$ by
$$(g \cdot f)(x) = g \cdot f (g^{-1} x), \, \forall g \in G$$
%

\begin{proof}

The operad substitution for $X_+ \sma \sC$ is given by
$$\xymatrix{
(X_+ \sma \sC(k)) \sma (X_+ \sma \sC(n_1)) \sma \dots \sma (X_+ \sma \sC(n_k)) \ar[d]^T \\
(X \times X^{\times k})_+ \sma (\sC(k) \sma \sC(n_1) \sma \dots \sma \sC(n_k)) \ar[d]^{p_+ \sma \gamma} \\
X_+ \sma \sC(\sum_i n_i)
}$$
Here $T$ is the permutation that collects the copies of $X_+$ together in the order that they appear, and $p: X \times X^{\times k} \to X$ projects away all but the first factor.  Note that unless $X$ is a point, the map $p$ keeps $X_+ \sma \sC$ from being a unital operad.

The operad substitution for $F(X_+, \; \sC)$ is
$$\xymatrix{
F(X_+, \; \sC(k)) \sma F(X_+, \; \sC(n_1)) \sma \dots \sma F(X_+, \; \sC(n_k)) \ar[d]^{smash} \\
F((X^{\times k+1})_+, \; \sC(k) \sma \sC(n_1) \sma \dots \sma \sC(n_k)) \ar[d]^{\Delta_+^*} \\
F(X_+, \; \sC(k) \sma \sC(n_1) \sma \dots \sma \sC(n_k)) \ar[d]^{\gamma_*} \\
F(X_+, \; \sC(\sum_i n_i))
}$$
Here $smash$ is the smash product of functions, and $\Delta: X \to X^{\times k+1}$ is the diagonal.

Both of these constructions are induced by lax monoidal functors, namely the smash product with $X_+$ and the function spectrum from $X_+$.  Since $\Delta$ makes $X$ into a \emph{cocommutative} co-monoid, the function spectrum functor $F(X_+, \cdot)$ is symmetric, so by Proposition \ref{lax_prop}, $F(X_+, \sC)$ is an operad.  Although $X_+ \sma \cdot$ is not symmetric, one can see directly that the above substitution map makes $X_+ \sma \sC$ an operad.


If $X$ and $\sC$ have $G$-actions, then the fact that the substitution for $X_+ \sma \sC$ is $G$-equivariant follows from the fact that both $\gamma$ and $p$ are.  Similarly, since $\Delta$ and $\gamma$ are equivariant, the substitution for $F(X_+, \; \sC)$ is equivariant. 

\end{proof}

There is an analogous statement for algebras: if $E$ is a $\sC$-algebra, then $X_+ \sma E$ is an $X_+ \sma \sC$-algebra, and $F(X_+, \; E)$ is a $F(X_+, \; \sC)$-agebra.  The proof is nearly identical. \\

\noindent {\it Proof of Corollary \ref{homotopy_corollary}.}  To show that $\sC_{bG}$ and $\sC^{hG}$ are operads (Corollary \ref{homotopy_corollary}) we examine the definition of each:
$$\sC_{bG} = (i^* i_* (EG_+ \sma \sC))^G$$
and
$$\sC^{hG} = F(EG_+, \; \sC)^G$$
The ``push-forward'' change of universe functor $i_*$ is interpreted here to be Elmendorf and May's replacement $I_{U^G}^U$ defined in \cite{elmay} (see also \cite{manmay}).  This functor is not just lax monoidal; in fact it is part of a strong symmetric monoidal equivalence of categories.  Consequently it carries naive $G$-operads to stable $G$-operads.  The functor $i^*$ is the right adjoint to $i_*$; it is lax monoidal by virtue of the general fact that the right adjoint to any strong symmetric monoidal functor is lax monoidal.

Therefore $\sC_{bG}$ and $\sC^{hG}$ are both formed by the composite of operations that produce operads.  The statement for algebras follows similarly. \\ \qed \\

The functoriality of these constructions deserves comment.  All constructions are covariantly functorial for maps of $G$-operads $\sC \to \sD$.  For the most part, they are all contravariantly functorial in the group $G$.

More specifically, the fixed point functor is contravariant for group homomorphisms $\phi: H \to G$; this extends to a natural transformation 
$$\phi^*: \sC^G \to \sC^H$$
between the $G$ and $H$ fixed point functors for naive $G$-operads $\sC$.

Furthermore, if we choose to define $EG$ as the two-sided bar construction $B(*, G, G)$, then the functor $G \mapsto EG$ is covariant in $G$; i.e., $\phi$ induces a map $E\phi: EH \to EG$.  Consequently the functor $\sC \mapsto F(EG_+, \sC)$ is contravariant in $G$.  Combining these facts, we see that
$$\sC \mapsto \sC^{hG}$$
is a contravariant functor in $G$.

%
%
%
%
%
%
%
%

For naive $G$-operads $\Sigma^{\infty} \sC_+$ which are suspension spectra of $G$-operads $\sC$, the naive construction of the transfer operads is, up to homotopy, a contravariant functor of subgroups of $G$.  The relative transfer 
$$\tau_H^G: \Sigma^{\infty} \sC(k)_{hG}^{Ad_G} \to \Sigma^{\infty} \sC(k)_{hH}^{Ad_H}$$
gives a map of operads in the stable homotopy category.  We expect that, using appropriate change of group and universe functors, one can show that $\sC_{bG}$ is covariantly functorial in $G$ for general naive $G$-operads $\sC$.

Identical remarks for algebras over these operads hold.  Further, in Lemma \ref{op_map_lemma} below, we will see that the norm map is a map of operads, allowing us to compare the transfer and homotopy fixed point operads.  The discussion above then applies to show that this comparison is coherent across subgroups of $G$.

\subsection{Fixed point operads in spaces} \label{space_fixed_op_section}

It is worth exploring a version of fixed point operads in the category of spaces, without first moving to the stable category.

\begin{proposition} \label{space_op_prop}

If $\sC$ is a $G$-operad, the fixed point spaces for the group action, $\sC(k)^G$ assemble into an operad $\sC^G$ in the category of topological spaces.

\end{proposition}

\begin{proof}

The equivariance of the substitution map $\gamma$ for $\sC$ gives a restriction
$$\gamma^{G}: (\sC(k) \times \sC(n_1) \times \dots \times \sC(n_k) )^{G} \to \sC(\sum_i n_i)^{G}$$
Moreover, it is clear that an element $(c_0, c_1, \dots, c_k)$ of the left hand side is fixed by the diagonal $G$ action if and only if each $c_i$ is fixed by the individual action of $G$ on the appropriate term of the operad.  Therefore
$$(\sC(k) \times \sC(n_1) \times \dots \times \sC(n_k) )^{G} = \sC(k)^G \times \sC(n_1)^G \times \dots \times \sC(n_k)^G$$
and $\gamma^G$ can be used as a substitution map for $\sC^G$.  That $\sC^G$ satisfies all the requisite axioms of an operad is a corollary of the fact that $\sC$ does. 

\end{proof}

One can similarly show that if $X$ is a $\sC \rtimes G$ algebra, $X^G$ is a $\sC^G$-algebra.

\begin{example}{\bf Little disks operads}

For the $k$-dimensional little disks operad $\sC_k$ with an action of $SO(k)$, the fixed points
$$\sC_k(n)^{SO(k)} = \emptyset$$
if $n>1$.  So $\sC_k^{SO(k)}$ is a fairly trivial operad.  However, taking closed subgroups of $SO(k)$ will give nontrivial results; since for instance $(\R^{i+j})^{SO(i)} = \R^j$ when $i>1$, there is an isomorphism of operads
$$\sC_{i+j}^{SO(i)} \cong \sC_j$$

\end{example}

\begin{example}{\bf Complex conjugation}

Alternatively, one may endow $\sC_{2k}$ with an action of $\Z / 2$, given by conjugation in $\R^{2k} = \C^k$.  Then
$$\sC_{2k}^{\Z / 2} \cong \sC_k$$

\end{example}

\begin{example}{\bf Operads of moduli spaces}

This second example suggests the following: consider the operad $\Mbar_\C$ whose $k^{\rm th}$ term $\Mbar_\C(k) = \overline{\sM}_{0, k+1}$ is the Deligne-Mumford compactification of the moduli space of $k+1$ marked points on the Riemann sphere $\C P^1$ (studied in \cite{km,gk,getz0}, amongst many others).  The substitution in the operad is given by attaching Riemann spheres at the marked points (giving a nodal curve).  $\Mbar_\C$ admits an action of $\Z / 2$ by complex conjugation, and the fixed point operad $\Mbar_\C^{\Z / 2} =: \Mbar_{\R}$ is called the \emph{mosaic operad} (studied in \cite{devadoss,ehkr} and others), with $k^{\rm th}$ term the real points of the moduli space:
$$\Mbar_{\R}(k) = \overline{\sM}_{0, k+1}(\R)$$

\end{example}

%
%

%
%
%
%
%
%

\section{Tate cohomology} \label{tate_section}

In this section we study the relation between the transfer and homotopy fixed-point operads and prove Proposition \ref{tate_proposition}.  Let $X$ be a naive $G$-spectrum.  For our purposes (as in \cite{gm}, section I.5), the norm map $n^G: X_{bG} \to X^{hG}$ is the composite
$$\xymatrix{
(i^* i_*(EG_+ \sma X))^G \ar[r]^-{q_*} & (i^* i_* X)^G = F(S^0, i^* i_* X)^G \ar[r]^-{p^*} & F(EG_+, \; i^* i_* X)^G
}$$
where $p:EG_+ \to S^0$ collapses $EG$ to a point, and $q = p \sma 1: EG_+ \sma X \to X$ projects $EG$ away.

As written, the target of $n^G$ is not $X^{hG}$, but rather $(i^* i_* X)^{hG}$.  However, the unit  
$$\eta: X \to i^* i_* X$$
of the $(i_*, i^*)$ adjunction is an equivariant map which is a nonequivariant equivalence (\cite{gm}, Lemma 0.1).  Such maps induce equivalences of homotopy fixed point spectra (noted in \cite{acd}).

\begin{remark}

Often (as in \cite{acd,klein_dual}) the norm map is defined as a map
$$EG_+ \sma_G (S^{Ad_G} \sma X) \to X^{hG}$$
Precomposing $n^G$ defined above with the equivalence $EG_+ \sma_G (S^{Ad_G} \sma X) \simeq X_{bG}$ gives precisely this definition.

\end{remark}

\begin{lemma} \label{op_map_lemma}

The norm map $n^G: \sC_{bG} \to \sC^{hG}$ is a map of operads.

\end{lemma}

\begin{proof}

We note that every spectrum in the composite
$$\xymatrix{
(i^* i_*(EG_+ \sma \sC(k)))^G \ar[r]^-{q_*} & (i^* i_* \sC(k))^G \ar[r]^-{p^*} & F(EG_+, \; i^* i_* \sC(k))^G
}$$
is the $k^{\rm th}$ term of an operad.  So to prove the lemma, we may show that $p^*$ and $q_*$ are maps of operads.  This follows from the functoriality of the construction in Lemma \ref{X_lemma}.  That is, for a map of spaces $f: X \to Y$, and an operad $\sD$, there are operad maps
$$\xymatrix{
X_+ \sma \sD \ar[r]^-{f_*} & Y_+ \sma \sD; & & F(X_+, \; \sD) & F(Y_+, \; \sD) \ar[l]_-{f^*}
}$$

Finally, we note that the unit of the $(i_*, i^*)$ adjunction gives a map of operads
$$\eta: \sC \to i^* i_* \sC.$$
Thus the equivalence $\sC^{hG} \simeq (i^* i_* \sC)^{hG}$ is given by a map of operads.  So $(i^* i_* \sC)^{hG}$ may be replaced with $\sC^{hG}$. 

\end{proof}

Recall that for a naive $G$-spectrum $X$, the generalized Tate cohomology $X^{tG}$ of $G$ with coefficients in $X$ is defined to be either the cofiber \cite{klein_dual,gm} or fiber \cite{acd} of the norm map
$$n^G: X_{bG} \to X^{hG}$$
Using either definition (since they differ only by a suspension), Proposition \ref{tate_proposition} follows from Lemma \ref{op_map_lemma} and the long exact sequence in stable homotopy groups for a (co)fibration sequence of spectra.

\section{Homotopy fixed points and homology} \label{dual_section}

We would like to compute the homology of the transfer and homotopy fixed point operads constructed in the previous sections.  The homology of $\sC_{bG}(k)$ is simply a suspension of the $G$-Borel equivariant homology of $\sC(k)$ (as in Corollary \ref{homology_corollary}).  The homology of $\sC(k)^{hG}$ is less straightforward, except when the Tate cohomology $\sC(k)^{tG}$ vanishes.  In this section we study more computable versions of the homology of homotopy fixed point spectra.

\subsection{Continuous homology of homotopy fixed point spectra} \label{pro_defn_section}

Let $h_*$ be a ring homology theory.  Recall that if $h$ is the spectrum representing $h_*$ and $X$ is any spectrum, then the $h_*$-homology of $X$ is defined by
$$h_*(X) := \pi_*(h \sma X)$$

\begin{definition}

For a ring homology theory $h_*$ and a naive $G$-spectrum $X$, the \emph{continuous homology of $X^{hG}$}, $h_*^{c}(X^{hG})$ is defined to be
$$h_*^{c}(X^{hG}) := \pi_*((h \sma X)^{hG})$$

\end{definition}

By definition, $h_*(X^{hG}) := \pi_*(h \sma (X^{hG}))$, differing somewhat from $h_*^{c}(X^{hG})$.  There is a natural map $e: h_*(X^{hG}) \to h_*^{c}(X^{hG})$ induced in $\pi_*$ by the composite
$$\xymatrix{
h \sma X^{hG} \ar[r]^-{\eta \sma 1} & h^{hG} \sma X^{hG} \ar[r]^-{i} & (h \sma X)^{hG}
}$$
Here $i$ is the lax monoidal natural transformation for the homotopy fixed point functor.  Since $h$ is a trivial $G$-spectrum, 
$$h^{hG} = F(EG_+, h)^G = F(BG_+, h)$$
and $\eta$ is induced by the projection $BG_+ \to S^0$.

%
%

The map $e$ is rarely an isomorphism.  Consider $X = S^0$:

\begin{proposition}

There is an isomorphism
$$h_*^{c}((S^0)^{hG}) \cong h^{-*}(BG)$$

\end{proposition}

\begin{proof}

The continuous homology $h_*^{c}((S^0)^{hG})$ is given by the homotopy groups of the spectrum
$$(h \sma S^0)^{hG} = F(EG_+, h)^G \simeq F(BG_+, h)$$
The last equivalence follows since the action of $G$ on $h$ is trivial.  But the cohomology of $BG$ is given by the homotopy groups of $F(BG_+, h)$. 

\end{proof}

In contrast,
$$h_*((S^0)^{hG}) = h_*(F(BG_+, S^0))$$
If $BG$ is a finite dimensional complex, then Spanier-Whitehead duality ensures that these two homologies are isomorphic.  However, when $BG$ is infinite dimensional (as is the case, for instance, when $G$ is compact Lie), one does not expect the homology of the Spanier-Whitehead dual $F(BG_+, S^0)$ to coincide with the (negative) cohomology of $BG$.

In the opposite extreme, the two homologies agree when $X$ is a free naive $G$-spectrum:

\begin{proposition} \label{pro_prop}

If $G$ is a compact Lie group and $X$ is a free naive $G$-CW spectrum, then $h_*(X^{hG})$ and $h_*^{c}(X^{hG})$ are mutually isomorphic to $h_*(X_{bG})$.

\end{proposition}

\begin{proof}

If the action of $G$ on $X$ is free, then it is also free on $h \sma X$.  Thus the Tate cohomology $(h \sma X)^{tG}$ vanishes, so that
$$(h \sma X)^{hG} \simeq (h \sma X)_{bG} \simeq EG_+ \sma_G (h \sma X \sma S^{Ad_G})$$
The action of $G$ on $h$ is trivial, so this last term is equivalent to
$$h \sma (EG_+ \sma_G (X \sma S^{Ad_G})) \simeq  h \sma X_{bG}$$
Taking homotopy groups then gives $h_*^{c}(X^{hG}) \cong h_*(X_{bG})$.  But freeness also implies that $X_{bG} \simeq X^{hG}$. 

\end{proof}

\subsection{Continuous homology of $LM^{hS^1}$} \label{pro_loop_section}

Examine the homotopy fixed point spectrum $LM^{hS^1}$.  By a result of Carlsson's \cite{carlssonloop}, if $M$ is simply connected and one completes at a prime $p$, there is a splitting
$$LM^{hS^1} \simeq \Sigma^{\infty} X_+ \lor \bigvee_{i=1}^{\infty} \Sigma \Sigma^{\infty} {LM_{hS^1}}_+$$
and the wedge is over $p$-adic valuations of the positive integers.  The norm map 
$$n^{S^1}: LM_{bS^1} \to LM^{hS^1}$$
is an equivalence of $LM_{bS^1} \simeq \Sigma \Sigma^{\infty} {LM_{hS^1}}_+$ with one of these wedge factors.

The (mod $p$) homology of $LM^{hS^1}$ is therefore quite complicated.  Notice, for instance, that if one takes $M$ to be a point, we see that after $p$-completion, 
$$(S^0)^{hS^1} \simeq S^0 \lor \bigvee_{i=1}^{\infty} \Sigma \Sigma^{\infty} BS^1_+$$
whose homology is substantially different from the negative cohomology of $BS^1$.

For naive $S^1$-spectra $X$, there is a homotopy fixed point spectral sequence
$$H^{-s}_{gp}(S^1, \pi_t(X)) \implies \pi_{s+t}(X^{hS^1})$$
(see, e.g., \cite{bruner_rognes}, remark 2.2b).  Taking $X$ to be $Y \sma h$ gives a conditionally convergent spectral sequence
$$H^{-s}_{gp}(S^1, h_t(Y)) \implies h^{c}_{s+t}(Y^{hS^1})$$
When $h = H\F_p$, this agrees with the homological homotopy fixed point spectral sequence which converges to the continuous homology $H^c_*(Y^{hS^1}, \F_p)$ (studied by Bruner-Rognes in \cite{bruner_rognes}).  In particular this gives a spectral sequence for computing the continuous homology of $LM^{hS^1}$:
$$Ext_{H_*(S^1)}^{-s}(\Z, H_t(LM)) \implies H_{s+t}^{c}(LM^{hS^1})$$
For examples of calculations of this nature using an analogous spectral sequence converging to $H_*(LM_{hS^1})$, we refer the reader to \cite{wes_string}.

\subsection{(Continuous) homology of operads} \label{pro_operad_section}

It is natural to study operads $\sC$ of topological spaces or spectra via their homology.  If $h_*$ is a ring homology theory, then $h_*(\sC) = \{ h_*(\sC(k)) \}$ is an operad in the category of $h_*$-modules.  We can see this as follows: if we assume that $h$ is a commutative $S$-algebra, then the functor
$$X \mapsto h \sma X$$
is lax symmetric monoidal; the natural transformation 
$$(h \sma X) \sma (h \sma Y) \to h \sma (X \sma Y)$$
is given by multiplication in $h$.  Furthermore, $\pi_*$ is lax symmetric monoidal, the natural transformation $\pi_*(X) \otimes \pi_*(Y) \to \pi_*(X \sma Y)$ being given by the smash product of maps.  Hence $h_*$ is a composite of lax symmetric monoidal functors, and thus lax symmetric monoidal; therefore $h_*(\sC)$ is an operad.  Notice that after applying $\pi_*$, $h$ need only be associative and commutative up to homotopy.

\begin{definition}

If $\sC$ is a naive $G$-operad, then the \emph{continuous homology of $\sC^{hG}$}, $h_*^{c}(\sC^{hG})$ is defined to be the operad (in the category of $h_*$-modules) whose $k^{\rm th}$ term is
$$h_*^{c}(\sC^{hG}(k)) = \pi_*((h \sma \sC(k))^{hG})$$

\end{definition}

The same arguments (in a different order) that prove that $h_*(\sC^{hG})$ is an operad serve to show that $h_*^{c}(\sC^{hG})$ is an operad.

\section{Examples of transfer and homotopy fixed point operads} \label{example_section}

In this section we study the transfer and homotopy fixed point operads for the $SO(2) = S^1$-action on $\sC_2$, and the $SU(2) = S^3$-action on $\sC_4$.  More specifically, we compute their associated homology operads which we show to be incarnations of the gravity operad.  We also study an operad $\unG$ constructed from products of a group $G$.  The associated transfer and homotopy fixed point operads appear to be related to the string topology of $BG$.  Since all of the naive $G$-operads considered in this paper will be suspension spectra of $G$-operads in spaces, we will frequently drop the $\Sigma^\infty \cdot_+$ from the notation.  

\subsection{The two-dimensional gravity operad} \label{2d_section}

We note that since $S^1$ is abelian, the adjoint action of $S^1$ on its Lie algebra is trivial, so 
$$
{\sC_2}_{bS^1}(k) \simeq ES^1_+ \sma_{S^1} (S^{Ad_{S^1}} \sma \Sigma^{\infty} \sC_2(k)_+) =  \Sigma \Sigma^{\infty} ({\sC_2(k)}_{hS^1})_+$$
We will denote the integral homology of this operad by ${e_2}_{S^1}$, echoing the convention $e_2 := H_*(\sC_2)$ in \cite{gj}; this is given by
$${e_2}_{S^1}(k) = \Sigma H_*({\sC_2(k)}_{hS^1})$$

As discussed in section \ref{pro_operad_section}, we will consider two possible interpretations of the homology of the homotopy fixed point operad:
$$\begin{array}{ccc}
H_*(\sC_2(k)^{hS^1}) & \;\;\; {\rm and} \;\;\; & H_*^{c}(\sC_2(k)^{hS^1})
\end{array}$$
Proposition \ref{pro_prop} and Lemma \ref{free_space_lemma} below will ensure that these homologies agree when $k>1$, though they will differ at $k=1$.  We will have more use for the continuous homology operad, so will define the operad $e_2^{S^1}$ as the continuous homology of $\sC_2^{hS^1}$:
$$e_2^{S^1}(k) = H_*^{c}(\sC_2(k)^{hS^1})$$

\begin{lemma} \label{free_space_lemma}

The action of $S^1$ on $\sC_2(k)$ is free when $k>1$.

\end{lemma}

This follows from the fact that the action of $S^1$ on $\R^2$ is free away from $0$, so is free on configurations of at least $2$ points.  There is a similar statement for the action of $SU(2)$ on $\sC_4(k)$.

When $X$ is a finite, free $G$-CW complex, the generalized Tate-cohomology $X^{tG}$ is trivial \cite{klein_dual} so this lemma, with Proposition \ref{tate_proposition}, has Corollary \ref{c2_corollary} as a corollary.  Additionally we have a homological analogue:

\begin{lemma} \label{free_lemma}

If $k>1$, $e_2(k) = H_*(\sC_2(k))$ is a free $H_*(S^1)$-module.

\end{lemma}

We will give a proof of this fact in section \ref{tech_section}.  Write the fundamental class of $S^1$ as $\Delta$, so that the Pontrjagin ring of $S^1$ is the exterior algebra
$$H_*(S^1) = \Lambda[\Delta]$$

\begin{corollary} \label{computation_corollary}

A computation of ${e_2}_{S^1}$ and $e_2^{S^1}$:

\begin{enumerate}

\item For $k>1$, ${e_2}_{S^1}(k) \cong e_2^{S^1}(k)$ is isomorphic to \label{ker_item}
$$\ker \Delta: H_*(\sC_2(k)) \to H_{*+1}(\sC_2(k)).$$

\item ${e_2}_{S^1}(1) \cong \Sigma H_*(BS^1)$.  It is a trivial ring. \label{item_2}

\item $e_2^{S^1}(1)$ is isomorphic to $H^{-*}(BS^1)$ as a ring. \label{item_3}

\end{enumerate}

\end{corollary}

\begin{proof}

To compute ${e_2}_{S^1}(k)$, we will use the naive description of the construction of the transfer operad and the Bousfield-Kan (or Borel) spectral sequence for the simplicial space $X_{hG}$: 
$$\Tor_*^{H_*(G)} (\Z, H_*(X)) \implies H_*(X_{hG})$$
For $X = \sC_2(k)$ ($k>1$) and $G = S^1$, Lemma \ref{free_lemma} implies that the spectral sequence collapses at the $E_2$-term, which is the quotient $H_*(\sC_2(k)) \otimes_{\Lambda[\Delta]} \Z$.  In homology, the transfer
$$\tau: \Sigma \Sigma^{\infty} {\sC_2(k)_{hS^1}}_+ \to \Sigma^{\infty} \sC_2(k)_+$$
carries $\Sigma H_*(\sC_2(k)) \otimes_{\Lambda[\Delta]} \Z$ isomorphically to 
$$\Delta(H_*(\sC_2(k))) = \ker \Delta$$
because, if $\pi: X \to X_{hS^1}$, then $\Delta = \tau_* \pi_*$, and $H_*(\sC_2(k)_{hS^1}) = \im \pi_*$.

By Lemma \ref{op_map_lemma}, the transfer is a map of operads, so away from the unary term this embeds ${e_2}_{S^1}$ as a suboperad of $e_2$.  An alternate proof of this fact follows from the computations in section \ref{tech_section}.

The unary term of any operad is a monoid; this induces the ring structure in homology.  Since $\sC_2(1) \simeq \R^2$ is $S^1$-equivariantly contractible to $\{ 0 \} \subseteq \R^2$, the computations of the homology in part \ref{item_2} follows.  That ${e_2}_{S^1}(1)$ is a trivial ring follows from the fact that it is concentrated in odd degrees.

The $S^1$-contractibility of $\sC_2(1)$ also implies that
$$\sC_2^{hS^1}(1) \simeq (S^0)^{hS^1} = F(BS^1_+, S^0)$$
whose continuous homology is the negative cohomology of $BS^1$.

Multiplication in $\sC_2^{hS^1}(1)$ is given by
$$
\xymatrix{ (\sC_2(1)_+)^{hS^1} \sma (\sC_2(1)_+)^{hS^1} \ar[r] \ar[d]_-{\simeq} & (\sC_2(1) \times \sC_2(1))_+^{hS^1} \ar[r]^-{\gamma^{hS^1}} \ar[d]_-{\simeq} & (\sC_2(1)_+)^{hS^1} \ar[d]_-{\simeq} \\
{S^0}^{hS^1} \sma {S^0}^{hS^1} \ar[r] & {S^0}^{hS^1} \ar[r]_-{\simeq} & {S^0}^{hS^1}}
$$
Under the identification of $(S^0)^{hS^1}$ with the Spanier-Whitehead dual of $BS^1_+$, the first map in the lower sequence is equivalent to the dual of the diagonal map.  Part \ref{item_3} follows. 

\end{proof}

We immediately obtain a proof of Proposition \ref{gravity_prop} via a comparison of part (\ref{ker_item}) and the discussion in \cite{getz2d} following Theorem 4.2.  An alternative proof is given in section \ref{tech_section}, where a direct equivalence
$$\sC_2(k)_{hS^1} \simeq \sM_{0, k+1} \leqno{(*)}$$
is given for $k>1$.  Interestingly, while there is not a good notion of the moduli space $\sM_{0, 2}$, one can define the moduli \emph{stack} $\sM_{0, 2}$ as the translation groupoid for the action of $PSL(2, \C)$ on the configuration space of $2$ points in $\C P^1$:
$$\sM_{0, 2} = [F(\C P^1, 2) /  PSL(2, \C)] \cong [\C / Aff(\C)]$$
Since $Aff(\C) = \C \rtimes \C^{\times}$ is homotopy equivalent to its subgroup $S^1$ and $\C$ is equivariantly contractible, the geometric realization of the stack $\sM_{0, 2}$ is equivalent to $BS^1$.  So in this sense, the transfer operad ${\sC_2}_{bS^1}$ allows us to extend Getzler's definition of the gravity operad in terms of $\sM_{0, k+1}$ to include $k=1$.

In any case, Proposition \ref{gravity_prop} allows us to refer to 
$$({e_2}_{S^1})_{>1} \cong (e_2^{S^1})_{>1}$$
as $Grav$.  Technically, Getzler's definition of $Grav$ is over $\C$, but it is easy to see that the description in \cite{getz2d} of $Grav(k)$ as $\ker \Delta$ could have been made integrally.

Examine the full operad $e_2^{S^1}$.  Let $A$ be an algebra over it; this endows $A$ with the structure of a Lie algebra from the embedding of operads $Lie \subseteq Grav \subseteq e_2^{S^1}$.  It also makes $A$ into a module over the ring
$$e_2^{S^1}(1) = H^{-*}(BS^1) = \Z [c]$$
where $c$ has dimension $-2$.

\begin{proposition}

Let $\theta \in e_2^{S^1}(k)$ for $k>1$, and $a_i \in A$ for $i = 1, \dots, k$.  Then
$$\theta(a_1, \dots, ca_i, \dots, a_k) = 0$$

\end{proposition}

\begin{proof}

Let $\theta' = \theta \circ_i c = \gamma(\theta; 1, \dots, c, \dots, 1)$; then
$$\theta(a_1, \dots, ca_i, \dots, a_k) = \theta'(a_1, \dots, a_k)$$
However, $e_2(k)$ is a free $H_*(S^1)$-module so $e_2^{S^1}(k)$ is a trivial $H^{-*}(BS^1)$-module.  Thus $\theta' = 0$. 

\end{proof}

\begin{corollary} \label{splitting_corollary}

An algebra $A$ over $e_2^{S^1}$, considered as a gravity algebra, contains a trivial gravity ideal $cA$.

\end{corollary}

Note that over a field, $A$ splits as the sum
$$A = cA \oplus A/cA$$
In particular we see that as a Lie algebra, $A$ splits as a $cA \oplus A/cA$, and that $cA$ is an abelian Lie algebra.  An analogous result holds for ${e_2}_{S^1}$; however, since ${e_2}_{S^1}(1)$ is already a trivial ring, it is not a noteworthy fact.

Consider the (usual, not continuous) homology operad of $\sC_2^{hS^1}$.  As mentioned above, $H_*(\sC_2^{hS^1}(k))$ agrees with $e_2^{S^1}(k)$ when $k>1$.  For $k=1$, using Carlsson's result \cite{carlssonloop}
$$H_*(\sC_2^{hS^1}(1), \F_p) = H_*(S^0 \lor \bigvee_{i=1}^{\infty} \Sigma \Sigma^{\infty} BS^1_+, \F_p)$$

Finally, we examine the fixed point operad $(i^* i_* \Sigma^{\infty} {\sC_2}_+)^{S^1}$.  Of course the fixed point operads $\sC_2^{S^1}$ (in spaces) and $(\Sigma^{\infty} {\sC_2}_+)^{S^1}$ are trivial by the freeness results described above.  We will therefore not consider them, and allow ourselves the notational simplification
$$\sC_2^{S^1} = (i^* i_* \Sigma^{\infty} {\sC_2}_+)^{S^1}.$$

A key fact is the tom Dieck splitting: for a $G$-space $X$, there is an equivalence
$$(i^* i_* \Sigma^{\infty} X_+)^G = (i^* \Sigma^{\infty}_G X)^G \simeq \bigvee_{(H) \leq G} \Sigma^{\infty} (X^H)_{hW(H)}^{Ad_{W(H)}}.$$
The wedge sum is over all conjugacy classes of closed subgroups $H \leq G$.  $W(H) = N(H)/H$ is the Weyl group of $H$ in $G$ which inherits a residual action on the fixed point space $X^H$.  In terms of this splitting, the projection map
$$p: X_{bG} = (i^* i_* \Sigma^{\infty} EG \times X_+)^G \to (i^* i_* \Sigma^{\infty} X_+)^G$$
carries $X_{bG}$ isomorphically onto the summand $\Sigma^{\infty} X_{hG}^{Ad_{G}}$ given by $H=0$.

\begin{proposition}

The map $p$ induces an equivalence $({\sC_2}_{bS^1})_{>1} \simeq (\sC_2^{S^1})_{>1}$.  Also, there is an equivalence
$$\sC_2(1)^{S^1} \simeq S^0 \lor \bigvee_{n \geq 0} \Sigma \Sigma^{\infty} B (S^1 / (\Z /n))_+$$

\end{proposition}

The group $S^1 / (\Z /n)$ is of course isomorphic to $S^1$, but this presentation reflects the tom Dieck splitting.

\begin{proof}

This proposition follows from the fact that, for $k>1$, $\sC_2(k)$ is a free $S^1$-space, and for $k=1$ it is $S^1$-equivariantly contractible. 

\end{proof}

\subsection{The four-dimensional gravity operad} \label{4d_section}

In this section we study the action of $SU(2)$ on the four-dimensional little disks operad $\sC_4$.  Here, $\sC_4$ becomes an $SU(2)$-operad via the representation $SU(2) \to SO(4)$.  One can define a \emph{four-dimensional gravity operad} $Grav^4$ of graded groups with an identical presentation as $Grav$, where the generators are taken to be of dimension $3$ instead of $1$.

\begin{proposition} \label{grav4_proposition}

There is an isomorphism $(H_*({\sC_4}_{bSU(2)}))_{>1} \cong Grav^4$.

\end{proposition}

We will reduce this computation to the two-dimensional computations in the previous section.  We need three key facts:

\begin{enumerate}

\item There is a ungraded ring isomorphism $H_*(S^1) \cong H_*(SU(2))$ (degrees are multiplied by $3$). \label{fact1}

\item There is a ungraded isomorphism $H_*(\sC_2(k)) \cong H_*(\sC_d(k))$ (degrees are multiplied by $d-1$). \label{fact2}

\item When $d=4$ in the previous, the isomorphism is one of $H_*(S^1) \cong H_*(SU(2))$-modules when one takes into account the isomorphism of fact \ref{fact1} above. \label{fact3}

\end{enumerate}

The first fact is straightforward, as both rings are exterior algebras on one generator.  The second fact follows from F. Cohen's computation of the homology of configuration spaces \cite{clm}.  For the third, recall that for $k>1$, Cohen's computations show that the cohomology of $\sC_d(k)$ is generated by pullbacks of the top-dimensional classes under the maps
$$\begin{array}{ccc}
\pi_{i, j}: \sC_d(k) \to \sC_d(2) \simeq S^{d-1}, & & i \neq j
\end{array}$$
which forget all but the $i^{\rm th}$ and $j^{\rm th}$ little disks.  These maps are $SO(d)$-equivariant, so fact \ref{fact3} reduces to the case $k=2$.  But there are equivariant equivalences 
$$\xymatrix{
\sC_2(2) \simeq_{S^1} S^1 & & \sC_4(2) \simeq_{SU(2)} SU(2)
}$$
so fact \ref{fact3} for $k=2$ is equivalent to fact \ref{fact1}. \\

\noindent {\it Proof of Proposition \ref{grav4_proposition}.} We compute $H_*({\sC_4(k)}_{bSU(2)})$ as we did in the two-dimensional case, using the Bousfield-Kan spectral sequence
$$\Tor_*^{H_*(SU(2))} (\Z, H_*(\sC_4(k))) \implies H_*({\sC_4(k)}_{hSU(2)}).$$
But thanks to the facts above, this is precisely the same computation that we made for $H_*({\sC_2(k)}_{hS^1})$.  Consequently there is a ungraded isomorphism of operads
$$(H_*({\sC_2}_{bS^1}))_{>1} \cong (H_*({\sC_4}_{bSU(2)}))_{>1}$$
in which degrees are scaled by a factor of $3$. \\ \qed

\subsection{Operads from loop spaces of classifying spaces} \label{loop_subsection}

In this section we offer an example of equivariant operads of a decidedly different flavor from those presented in the previous sections.

\begin{definition}

For a group $G$, define an operad $\unG$ by $\unG(k):= G^{\times k}$, with substitution map 
$$\gamma: \unG(k) \times \unG(n_1) \times \dots \times \unG(n_k) \to \unG(\sum n_i)$$
defined by
$$\gamma(g_1, \dots, g_k; h_1^1, \dots, h_{n_1}^1, \dots, h_1^k, \dots, h_{n_k}^k) = (g_1 h_1^1, \dots, g_1 h_{n_1}^1, \dots, g_k h_1^k, \dots, g_k h_{n_k}^k)$$

\end{definition}
Equivalently, $\unG$ is the semidirect product $Comm \rtimes G$, where $Comm$ is the commutative operad, given the trivial $G$-action.

Define an action of $G$ on $\unG(k)$ by (diagonal) conjugation on each factor.  As $\gamma$ is equivariant with respect to this action, it makes $\unG$ into a $G$-operad.  We are thus entitled to form the transfer and homotopy fixed point operads, $\unG_{bG}$ and $\unG^{hG}$.

Recall that when $G$ acts on itself by conjugation, there is a (fibrewise) equivalence between the Borel construction and the free loop space of the classifying space of $G$:
$$\xymatrix{
G \times_G EG \ar[r]^-{\simeq} \ar[d]_-{proj_2} & LBG \ar[d]^-{ev} \\
BG \ar[r]^-{=} & BG
}$$
The following is a simple application of this fact:

\begin{proposition} \label{LBG_op_prop}

For a compact Lie group $G$, there is an equivalence
$$\unG_{bG}(k) \simeq \Sigma^{\infty} (LBG \times_{BG} LBG \times_{BG} \dots \times_{BG} LBG)^{Ad_G}$$
where there are $k$ factors in the iterated fiber product.

\end{proposition}

Recall that the unary term of an operad is always a monoid; we see then that $\Sigma^{\infty} LBG^{Ad_G}$ obtains the structure of a ring spectrum.  So if $Ad_G$ is an orientable vector bundle, then $\Sigma^{\dim G} H_*(LBG)$ forms a ring.  

In \cite{acg} Abbaspour, Cohen, and Gruher defined string topology operations (in particular the loop product) on $H_*(LBG)$ when $G$ is a Poincar\'e duality group, and described the operation in group theoretic terms.  Gruher and Salvatore extended this to compact Lie groups $G$ in \cite{grusalv} using geometric methods which rely upon the fact that $BG$ is a colimit of finite dimensional manifolds.  

These constructions produce markedly different ring structures on $H_*(LBG)$ from the one defined by Proposition \ref{LBG_op_prop}, as is apparent from a simple examination of the degree shifts involved.  Perhaps a more compelling argument is as follows: Recall (as in Corollary \ref{computation_corollary}) that $\Sigma^{\infty} BG^{Ad_G}$ is a ring spectrum via the diagonal transfer (for finite groups this structure was studied in detail in the $K(n)$-local category by Strickland in \cite{strick}).

\begin{proposition} \label{abelian_prop}

If $G$ is abelian, then there is an equivalence of ring spectra
$$\Sigma^{\infty} G_+ \sma \Sigma^{\infty} BG^{Ad_G} \simeq \unG_{bG}(1) \simeq \Sigma^{\infty} LBG^{Ad_G}$$
where the left hand side is a (smash) product of ring spectra.

\end{proposition}

\begin{proof}

Since $G$ is abelian, $G \times_G EG = G \times BG$.  We will use the naive description of the operadic structure on $\unG_{bG}$ (and hence the monoidal structure on $\Sigma^{\infty} LBG^{Ad_G}$).  In that description, the multiplication in $\Sigma^{\infty} LBG^{Ad_G}$ is given by
$$\xymatrix{
(G_+ \sma BG^{Ad_G}) \sma (G_+ \sma BG^{Ad_G}) \ar[r]^-{\tau_{\Delta(G)}^{G \times G}} & (G \times G)_+ \sma BG^{Ad_G} \ar[r]^-{\mu \sma 1} & G_+ \sma BG^{Ad_G}
}$$
(we drop the $\Sigma^{\infty}$ for brevity).  Because the action of $G$ on itself by conjugation is trivial, this is precisely the smash product of the multiplication in $G_+$ with the multiplication in $\Sigma^{\infty} BG^{Ad_G}$. 

\end{proof}

If one takes $G=S^1 =: \T$, then, as we have seen above, $H_*({B\T}^{Ad_{\T}})$ is a trivial ring, so the ring structure on $H_*({LB\T}^{Ad_{\T}})$ is a tensor product of a trivial ring with an exterior algebra on one generator of dimension $1$.  In constrast, it is shown in \cite{grusalv} that Gruher and Salvatore's loop product for ${LB\T}^{-TB\T}$ produces the pro-ring
$$H_*({LB\T}^{-TB\T}) \cong \Lambda(t) \otimes \Z[[c]]$$
where $t$ has dimension $1$, and $c$ has dimension $-2$.

Notice, however, that the same sort of arguments in Proposition \ref{abelian_prop} give the following fact:

\begin{proposition}

If $G$ is abelian, then there is an equivalence of ring spectra
$$\Sigma^{\infty} G_+ \sma F(\Sigma^{\infty} BG_+, S^0) \simeq \unG^{hG}(1)$$
where the left hand side is a smash product of ring spectra.

\end{proposition}

\noindent Here, the Spanier-Whitehead dual $F(\Sigma^{\infty} BG_+, S^0)$ is given the usual ring structure which induces the cup product in $H^*BG$.

When $G = \T$ we see that the continuous homology $H_*^{c}(\unG^{hG}(1))$ produces the same ring as Gruher and Salvatore's construction.  In future work \cite{gw} we will show that this is in fact true for all compact Lie groups $G$.

We can also consider fixed point operads.  In the category of spaces, the fixed points of the conjugation action of $G$ on itself is the center of $G$:
$$G^G = Z(G)$$
This extends to products: $(G^{\times k})^G = Z(G)^{\times k}$.  Consequently, in the category of spaces, there is an isomorphism of operads
$$\unG^G = \underline{Z(G)}$$

In the stable category, we can consider the operad
$$(i^* \Sigma^{\infty}_G \unG_+)^G = (i^* i_* \Sigma^{\infty} \unG_+)^G$$
For a subgroup $H \leq G$, the fixed point set $G^H = C(H)$ is the centralizer of $H$ in $G$.  So, via the tom Dieck splitting, we may determine the homotopy type of the terms of the operad:
$$(i^* \Sigma^{\infty}_G \unG_+)^G(k) \simeq \bigvee_{(H) \leq G} \Sigma^{\infty} (C(H)^{\times k})_{hW(H)}^{Ad_{W(H)}}$$
where the sum is over conjugacy classes of closed subgroups $H \leq G$.

\section{String homology} \label{string_homology_section}

In this section, let $M$ be a closed $d$-manifold, and let $h_*$ be a ring homology theory with respect to which $M$ is orientable.  The Cohen-Jones construction \cite{cj} gives an action of the cactus operad $Cacti$ on $LM$ in the category of correspondences of spaces (see also \cite{cv}).  This is not quite sufficient to make $LM$ (or even $LM^{-TM}$) an algebra over $Cacti$ in spaces (or spectra).  Nonetheless, their construction may be used to give $h_*(LM)$ the structure of an $h_*(Cacti)$-operad.  

Since $Cacti$ is homotopy equivalent to the framed disk operad $\sC_2 \rtimes S^1$, we would like to adapt Cohen and Jones' framework along the equivariant lines described above to give $h_*(LM_{bS^1})$ the structure of an algebra over the transfer operad $h_*({\sC_2}_{bS^1})$.

However, the general machinery set up in the previous sections cannot be directly applied to string homology.  Two substantial modifications are needed to circumvent certain technical issues.  The first problem is the fact that the $Cacti$ action is through correspondences.  As in \cite{cj} this can be dealt with by an appropriate \emph{umkehr} map.  The second problem is that while $Cacti \simeq \sC_2 \rtimes S^1$, $Cacti$ is \emph{not} a semidirect product of an operad and a group.  Rather, as Kaufmann \cite{kaufcacti} has pointed out, it is a \emph{bi-crossed product} of a suboperad $Cact$ with $S^1$.  

In sections \ref{cactus_subsection} and \ref{cactus2_section} we explore the impact of this subtlety, defining an operad $Cact_{bS^1}$ in analogy with the transfer operad ${\sC_2}_{bS^1}$.  In \ref{mod_subsection} we use $Cact_{bS^1}$ to give the operad action described in Theorem \ref{string_homology_theorem}.  In section \ref{bracket_subsection} we study this action when $h$ is taken to be singular homology, and identify the Chas-Sullivan string bracket in terms of the operad ${e_2}_{S^1}$.  We develop (homotopy) fixed point versions in parallel with this discussion.


\subsection{Equivariance of the cactus operad} \label{cactus_subsection}

We will not belabor the definition of the cactus operad $Cacti$; we refer the reader to several careful treatments of it \cite{voruniv,cv,kaufcacti}.  Recall that the $k^{\rm th}$ space $Cacti(k)$ consists of isotopy classes of $k$-lobed cacti (configurations of $k$ circles in the plane whose dual graph is a tree), along with (inner) marked points on each circle, and an outer marked point on the cactus itself.  A cactus $c \in Cacti(k)$ determines a pinching map $\nabla_c: S^1 \to c$.

The partial substitution map for $Cacti$ 
$$\circ_i: Cacti(m) \times Cacti(n) \to Cacti(m+n-1)$$
is given by replacing the $i^{\rm th}$ lobe of the first cactus with the second cactus, via the pinching map.

Kaufmann \cite{kaufcacti} defines a suboperad $Cact < Cacti$ of \emph{spineless cacti}.  One definition (equivalent to his) of $Cact(k)$ is as follows:

\begin{definition}

Let $Cact(k)$ be the subspace of $Cacti(k)$ obtained by requiring that the inner marked point on the $i^{\rm th}$ circle be the first point (with respect to a clockwise orientation on $S^1$) on that circle in the image of the pinching map.

\end{definition}

Thus the inner marked point on the circle containing the outer marked point is also the outer marked point, and every other inner marked point is at intersections of lobes of the cactus.  Throughout this section we will use the abbreviation $C = Cact$.

\begin{theorem} \label{cact_is_e2_thm}

\cite{kaufcacti} $Cact$ is an $E_2$ operad.

\end{theorem}

There is an action of $S^1$ on $C(k)$ given by rotating the outer marked point clockwise around the cactus.  Unlike $\sC_2$, this does not make $C$ an $S^1$-operad; the operad substitution
$$\gamma: C(k) \times C(n_1) \times \dots \times C(n_k) \to C(\sum n_i)$$
is \emph{not} a map of $S^1$-spaces when $S^1$ acts diagonally on the left.  Consequently, we cannot use Corollary \ref{homotopy_corollary} to give an operad structure on $C_{bS^1}$ and $C^{hS^1}$.  We can get around this difficulty.

Kaufmann defines the \emph{homotopy diagonal defined by a spineless cactus}; this can be thought of as a map
$$\Delta: C(k) \to Map(S^1, (S^1)^{\times k})$$
For a cactus $c$, $\Delta(c):S^1 \to (S^1)^{\times k}$ is defined as a piecewise map gotten from going around the cactus clockwise, starting at the outer marked point.  $\Delta(c)$ is constant in the $i^{\rm th}$ component until the pinching map reaches the $i^{\rm th}$ lobe; then $\Delta(c)$ travels around the $i^{\rm th}$ term of $(S^1)^{\times k}$ until the pinching map reaches the next lobe.  We refer the reader to Figure 9 of \cite{kaufcacti} for a descriptive picture of the image of $\Delta(c)$.


We reinterpret Kaufmann's definition as follows:

\begin{definition} \label{hd_def}

Define the \emph{homotopy diagonal action} of $S^1$ on 
$$C(k) \times C(n_1) \times \dots \times C(n_k)$$
as follows: for $\theta \in S^1$, $c \in C(k)$, and $c_i \in C(n_i)$, let
$$\theta \cdot (c; \; c_1, \dots, c_k) = (\theta \cdot c; \; \Delta(c)(\theta) \cdot (c_1, \dots, c_k))$$
where, on the righthand side, $\cdot$ indicates the $S^1$ action on cacti by rotating the outer marking.

%
%
%

Write $\Delta(c)(\theta) = (\Delta_1(c)(\theta), \dots, \Delta_k(c)(\theta))$.  For $k, l \geq 1$ and $1 \leq i \leq k$, define the \emph{$i^{th}$ homotopy diagonal action} of $S^1$ on $C(k) \times C(l)$ by
$$\theta \cdot (c, d) = (\theta \cdot c, (\Delta_i(c)(\theta)) \cdot d)$$
We note that the $i^{\rm th}$ homotopy diagonal action is the equal to the homotopy diagonal action in the case $n_i = l$, and $n_j = 1$ if $j \neq i$.

%

\end{definition}

\begin{lemma} \label{action_lemma}

The homotopy diagonal and $i^{th}$ homotopy diagonal are group actions for every $i$.  Moreover, all are homotopic (through group actions) to the diagonal action.

\end{lemma}

\begin{proposition} \label{equiv_prop}

With respect to the $i^{th}$ homotopy diagonal action on the domain, the $i^{th}$ substitution
$$\xymatrix{
C(k) \times C(l) \ar[r]^-{\circ_i} & C(k+l-1)
}$$
is $S^1$-equivariant.  
%

\end{proposition}

The proofs of Lemma \ref{action_lemma} and Proposition \ref{equiv_prop} are somewhat technical and are relegated to section \ref{tech_section}.  

\subsection{Transfer and homotopy fixed point cactus operads} \label{cactus2_section}

Using Proposition \ref{equiv_prop} we may construct analogues of the transfer and homotopy fixed point operads for the $S^1$ action on $C$.  

\begin{theorem} \label{equiv_cactus_theorem}

The two families of spectra
$$
\begin{array}{rcl}
C_{bS^1} & := & \{ \Sigma (\Sigma^{\infty} C(k)_+)_{hS^1} \} \\
C^{hS^1} & := & \{ (\Sigma^{\infty} C(k)_+)^{hS^1} \}
\end{array}
$$
define operads in the stable homotopy category.

\end{theorem}

While for our purposes it is sufficient to show that these are operads up to homotopy, it would be quite surprising to find that they cannot be made into strict operads.  Further, we expect that a strict version of Theorem \ref{same_at_last_thm} below also holds.

We will call $C_{bS^1}$ the \emph{cactus transfer operad} and $C^{hS^1}$ the \emph{homotopy fixed point cactus operad}.  As was the case for the little disks operad, one can see that $(C_{bS^1})_{>1}$ and $(C^{hS^1})_{>1}$ are equivalent using Tate cohomology.

\begin{remark}

Since the circle action on $C(k)$ is by rotating the outer marking, one can picturesquely describe the terms of $C_{bS^1}$ as spaces of cacti without any markings at all.  In this description, the substitution maps are the substitution maps of $C$, summed over all possible outer markings.

\end{remark}

\begin{proof}

We will prove the theorem only for $C_{bS^1}$; the reader can make the appropriate modifications to prove that $C^{hS^1}$ is an operad.  We give an operad structure by specifying the $i^{\rm th}$ substitution $\uncirc_i$ for $C_{bS^1}$.

For each $1 \leq i \leq k$, equip $C(k) \times C(l)$ with an action of $S^1 \times S^1$ where the first factor acts via the $i^{\rm th}$ homotopy diagonal action, and the second factor acts on $C(l)$; i.e.,
$$(\theta, \phi) \cdot (c, d) = (\theta \cdot c, (\phi \cdot \Delta_i(c)(\theta)) \cdot d) \leqno{(*)}$$
The fact that $S^1$ is abelian ensures that this is in fact a group action.

Recall from Lemma \ref{action_lemma} that the $i^{\rm th}$ homotopy diagonal action is homotopic to the diagonal action through group actions.  Consequently the action in $(*)$ is homotopic through group actions to the action
$$(\theta, \phi) \cdot (c, d) = (\theta \cdot c, (\phi \cdot \theta) \cdot d)$$
The homotopy orbits of this action can be identified with $C(k)_{hS^1} \times C(l)_{hS^1}$ by first taking orbits with respect to $\phi$, then $\theta$.  Therefore there is a homotopy equivalence
$$\xymatrix{
j: C(k)_{hS^1} \times C(l)_{hS^1} \ar[r]^-{\simeq} & (C(k) \times C(l))_{h(S^1\times S^1)} 
}$$
where the homotopy orbits in the target are with respect to the action $(*)$.  Restricting to the the action of $\Delta_i(S^1) = S^1 \times 1 \leq S^1 \times S^1$ which acts by the $i^{\rm th}$ homotopy diagonal action, we also have an equivalence
$$\xymatrix{
j': (C(k) \times C(l))_{h\Delta(S^1)} \ar[r]^-{\simeq} & (C(k) \times C(l))_{h\Delta_i(S^1)} 
}$$
Here, homotopy orbits taken with respect to $\Delta(S^1)$ are via the diagonal action.

Then one may define $\uncirc_i$ for $C_{bS^1}$ by the diagram
$$\xymatrix{
\Sigma C(k)_{hS^1} \sma \Sigma C(l)_{hS^1} \ar[r]^-{\Sigma^2 j} \ar@{-->}[d]_{\uncirc_i} & \Sigma^2 (C(k) \times C(l))_{h(S^1\times S^1)} \ar[d]^-{\tau} \\
\Sigma C(k+l-1)_{hS^1} & \Sigma (C(k) \times C(l))_{h\Delta_i(S^1)} \ar[l]^-{\Sigma (\circ_i)_{hS^1}} \\
}$$
For brevity, we are dropping $\Sigma^{\infty} \cdot_+$ from our notation.  The relative transfer 
$$\tau = \tau_{\Delta_i(S^1)}^{S^1 \times S^1}$$
is with respect to the factor which acts by the $i^{\rm th}$ homotopy diagonal action. The map $(\circ_i)_{hS^1}$ exists by virtue of Proposition \ref{equiv_prop}.

The spectrum $C_{bS^1}(k)$ becomes a naive $\Sigma_k$-spectrum since the action of $\Sigma_k$ on $C(k)$ commutes with the $S^1$-action.  To show that this $\Sigma_k$-action and collection of substitution maps make $C_{bS^1}$ into an operad up to homotopy, we must show that substitution is associative and behaves correctly with respect to the $\Sigma_k$-action.  This follows from the next lemma, since $\W_{bS^1}$ and $\W^{hS^1}$ are operads. 

\end{proof}

\begin{lemma}

There exists an $S^1$-operad $\W$ and $\Sigma_k$-equivariant homotopy equivalences
$$\xymatrix@1{
\W_{bS^1}(k) \ar[r]^-{\simeq} & C_{bS^1}(k) & and & \W^{hS^1}(k) \ar[r]^-{\simeq} & C^{hS^1}(k)
}$$
that preserve operadic substitution up to homotopy.

\end{lemma}

\begin{proof}

Theorem \ref{cact_is_e2_thm} (proven using Fiedorowicz's recognition principle \cite{fied}) asserts that $\sC_2$ and $C$ are connected by some chain of equivalences.  Alternatively, this is given by a single equivalence
$$\psi: W(\sC_2) \to C,$$
where $W$ is the cofibrant replacement functor in the category of operads \cite{bv}.  By general nonsense, $\W := W(\sC_2)$ is an $S^1$-operad.

Similarly, there is an equivalence of the framed little disks $\sC_2 \rtimes S^1$ with the full cactus operad $Cacti$; this is a theorem of Voronov's \cite{voruniv}.  We will employ Kaufmann's argument: using the equivalence of $\sC_2$ and $Cact$, one immediately gets an equivalence of semidirect products $\sC_2 \rtimes S^1 \simeq Cact \rtimes S^1$.  The semidirect product $Cact \rtimes S^1$ is, in Kaufmann's terms, a ``quasi-operad;'' it has a substitution map that is not associative.  However, the homotopy between the diagonal action and the homotopy diagonal action gives an equivalence of quasi-operads $Cact \rtimes S^1 \simeq Cact \bowtie S^1 =: Cacti$ (the $\bowtie$ symbol denotes ``bicrossed product'').  As such, we have the following homotopy equivalences of quasi-operads:
$$\xymatrix{
\W \rtimes S^1 \ar[r]^-{\psi \rtimes 1} & Cact \rtimes S^1 \ar[r]^-{\simeq} & Cact \bowtie S^1 = Cacti
}$$

Consider the following diagram.  To define the leftmost square, use the fact that $(\W \rtimes S^1) (1) \simeq S^1 \simeq (C \rtimes S^1)(1)$.
$$
\xymatrix{
S^1 \times \W(k) \ar[r]^-{\subseteq} \ar[d]_-{1 \times \psi} & (\W \rtimes S^1) (1) \times (\W \rtimes S^1) (k) \ar[r]^-{\gamma} \ar[d]_-{(\psi \rtimes 1) \times (\psi \rtimes 1)}& (\W \rtimes S^1) (k) \ar[d]^-{\psi \rtimes 1} \ar[r] &  \W(k) \ar[d]^-{\psi} \\
S^1 \times C(k) \ar[r]_-{\subseteq} & (C \rtimes S^1)(1) \times (C \rtimes S^1)(k) \ar[r]_-{\gamma} & (C \rtimes S^1)(k)  \ar[r] & C(k)
}$$
The middle square commutes since $\psi \rtimes 1$ is a map of quasi-operads.  The rightmost square commutes since both horizontal maps project the $(S^1)^{\times k}$ factor away.

By definition, the horizontal compositions induce the $S^1$ action on the $k^{\rm th}$ term of each operad, so we see that $\psi$ is an equivariant map.  Therefore $\psi$ induces an equivalence of spectra
$$\psi_{bS^1}: \W_{bS^1}(k) \to C_{bS^1}(k)$$
Since $\psi$ is a map of operads, it, and hence $\psi_{bS^1}$ are $\Sigma_k$-equivariant. Thus, the result will follow if we show that $\psi_{bS^1}$ preserves substitution; that is, this homotopy commutes:
%
$$
\xymatrix{
\Sigma \W(k)_{hS^1} \sma \Sigma \W(l)_{hS^1} \ar[dd]_-{\tau_\Delta} \ar[r]^-{(\psi_{bS^1})^{\sma 2}} & \Sigma C(k)_{hS^1} \sma \Sigma C(l)_{hS^1} \ar[dd]_-{\tau_\Delta} \ar[r]^-{=}  & \Sigma C(k)_{hS^1} \sma \Sigma C(l)_{hS^1} \ar[d]^-{j}  \\
 & & \Sigma^2 (C(k) \times C(l))_{hS^1\times S^1} \ar[d]^-{\tau} \\
\Sigma (\W(k) \times \W(l))_{h\Delta(S^1)} \ar[d]_-{\Sigma (\circ_i)_{hS^1}} \ar[r]^-{(\psi^{\times 2})_{bS^1}} &\Sigma (C(k) \times C(l))_{h\Delta(S^1)}  \ar[r]^-{j'} & \Sigma (C(k) \times C(l))_{h\Delta_i(S^1)} \ar[d]^-{\Sigma (\circ_i)_{hS^1}} \\
\Sigma \W(k+l-1)_{hS^1} \ar[rr]_-{\psi_{bS^1}} & & \Sigma C(k+l-1)_{hS^1} \\
}$$
%
We are using the ``naive'' substitution map on $\W_{bS^1} \simeq \Sigma \W_{hS^1}$ of Section \ref{sketch_subsection} since we are only concerned with structures up to homotopy.  The upper left square commutes since transfers homotopy commute with equivariant maps.  The upper right pentagon commutes since $j$ is induced by $j'$.  The lower pentagon commutes up to homotopy since $\psi$ and $\psi \rtimes 1$ are maps of quasi-operads.

We leave it to the reader to adapt this proof to an equivalence of the homotopy fixed point operads. 

\end{proof}

Since the operads ${\sC_2}_{bS^1}$ and ${W(\sC_2)}_{bS^1} = \W_{bS^1}$ are obviously equivalent (and likewise for homotopy fixed point operads), the following is a corollary to the previous lemma.

\begin{theorem} \label{same_at_last_thm}

In the stable homotopy category, there are equivalences of operads
$$\xymatrix{
{\sC_2}_{bS^1} \simeq Cact_{bS^1} & {\it and} & \sC_2^{hS^1} \simeq Cact^{hS^1}
}$$

\end{theorem}


\subsection{Applications to string homology} \label{mod_subsection}

Our goal in this section is to use the constructions in the previous sections to prove Theorem \ref{string_homology_theorem}.  In view of Theorem \ref{same_at_last_thm}, the reader may reasonably wonder why  we discussed the cactus transfer operad at all.  The reason is that the Cohen-Jones construction of the Batalin-Vilkovisky structure on $H_*(LM)$ goes through the cactus operad, and is not defined naturally for the little disk operad.

Recall that, for each $k>0$, Cohen and Jones define maps
$$\xymatrix{ Cact(k) \times LM^{\times k} & Cact(k)M \ar[r]^-{concat} \ar[l]_-{\supseteq} & LM }$$
Actually, they did this for $Cacti(k)$, not $Cact(k)$, but the construction restricts to that suboperad. Here $Cact(k)M$ is the subspace of $Cact(k) \times LM^{\times k}$ consisting of $k$-lobed cacti $c$ and loops $f_1, \dots, f_k \in LM$ that coincide at points of intersection of the cactus (and therefore define a map from $c$ into $M$).  The map $concat$ is the concatenation of the loops $f_1, \dots, f_k$ via the pinching map $\nabla_c: S^1 \to c$.

Let $S^1$ act on $Cact(k) \times LM^{\times k}$ via the homotopy diagonal action; for a cactus $c \in Cact(k)$, loops $f_1, \dots, f_k \in LM$, $\theta \in S^1$, let
$$\theta \cdot (c; f_1, \dots, f_k) = (\theta \cdot c; \Delta(c)(\theta) \cdot (f_1, \dots, f_k))$$
The proof of Lemma \ref{action_lemma} adapts immediately to show that this is in fact a group action.  Further, this action tautologically preserves the fact that the loops agree at intersection points prescribed by the cactus.  Thus the subspace $Cact(k)M$ inherits the action.  Let $\nu$ denote the normal bundle of $Cact(k)M \subseteq Cact(k) \times LM^{\times k}$.

\begin{lemma} \label{ptmap_lemma}

There is a Pontrjagin-Thom collapse map
$$\pi: (Cact(k) \times LM^{\times k})_{hS^1} \to (Cact(k)M_{hS^1})^{\nu}$$

\end{lemma}

\begin{proof}

First we note that the normal bundle $\nu$ is an $S^1$-equivariant bundle.  To see this, let $c \in Cact(k)$, and let $r$ be the number of intersections (with multiplicity) in $c$.  Recall that in the stratum $U \subseteq Cact(k)M$ of cacti homeomorphic to $c$, $\nu$ is a pullback of various copies of $TM$ along the map
$$ev: U \to M^{\times r}$$
given by evaluation of the loops at the points of intersection of the cactus.  As pointed out above, the action of $S^1$ on $Cact(k)M$ preserves the evaluation of loops at the intersection points of the cactus.  Thus $ev$ is $S^1$-invariant, so the pullback bundle is $S^1$-equivariant.

The subspace $Cact(k)M \subseteq Cact(k) \times LM^{\times k}$ is $S^1$-invariant, so the Pontrjagin-Thom collapse
$$Cact(k) \times LM^{\times k} \to Cact(k)M^{\nu}$$
is equivariant (being either the identity or constant), and therefore extends to
$$(Cact(k) \times LM^{\times k})_{hS^1} \to (Cact(k)M^{\nu})_{hS^1}$$
But since $\nu$ is an equivariant bundle, the target of the map may be written as $(Cact(k)M_{hS^1})^{\nu}$. 

\end{proof}

Extend the action of $S^1$ on $Cact(k) \times LM^{\times k}$ to an action of $(S^1)^{\times k+1}$ as in the construction of the cactus transfer operad; the first copy of $S^1$ acts via the homotopy diagonal map, and the remaining $k$ act termwise on $LM^{\times k}$.

\begin{definition}

Let $\tau$ be the relative transfer map for the subgroup $S^1 \times \{ 1 \}^{\times k} < (S^1)^{\times k+1}$ which acts by the homotopy diagonal action:
$$\tau: \Sigma^{k+1} \Sigma^{\infty} (Cact(k) \times LM^{\times k})_{h(S^1)^{\times k+1} \, +} \to \Sigma \Sigma^{\infty} (Cact(k) \times LM^{\times k})_{hS^1 \, +}$$
Also let $T$ denote the $h_*$-Thom isomorphism for $\nu$.  It is apparent that $concat$ is an $S^1$-equivariant map.  Thus one may give $\Sigma^{1-d} h_*(LM_{hS^1})$ an action of the cactus transfer operad by the map $\unrho$:
$$\unrho: \Sigma h_*(Cact(k)_{hS^1}) \otimes (\Sigma^{1-d} h_*(LM_{hS^1}))^{\otimes k} \to \Sigma^{1-d} h_*(LM_{hS^1})$$
which is the composite
$$\unrho = (concat_{hS^1})_* \circ T \circ \pi_* \circ \tau_*$$

\end{definition}

Here we are implicitly using the fact that
$$\Sigma^{k+1} \Sigma^{\infty} (Cact(k) \times LM^{\times k})_{h(S^1)^{\times k+1} \, +} \cong \Sigma \Sigma^{\infty} Cact(k)_{hS^1 \, +} \sma (\Sigma \Sigma^{\infty} LM_{hS^1 \, +})^{\sma k}$$
and the fact that, since $h_*$ is a ring homology theory, there is a map (though not necessarily an isomorphism)
$$h_*(X) \otimes h_*(Y) \to h_*(X \sma Y)$$

For comparison, we define an operad action of $h_*(Cact^{hS^1})$ on $\Sigma^{-d} h_*(LM^{hS^1})$  via a map
$$\ovrho: h_*(Cact(k)^{hS^1}) \otimes (\Sigma^{-d} h_*(LM^{hS^1}))^{\otimes k} \to \Sigma^{-d} h_*(LM^{hS^1})$$
by $\ovrho = (concat^{hS^1})_* \circ T \circ \pi_* \circ \phi_*$, where
$$\phi: (Cact(k) \times LM^{\times k})^{h(S^1)^{\times k+1}} \to (Cact(k) \times LM^{\times k})^{hS^1}$$
is the forgetful map to the homotopy fixed point spectrum for the homotopy diagonal action.

We will need the following technical lemma.  Its proof is straightforward, since $\pi$ is either the identity or constant.

\begin{lemma} \label{equivpt_lemma}

Let $N$ be a $G$-space and let $M \subseteq N$ be a $G$-subspace of finite codimension with $G$-equivariant normal bundle $\nu$.  For a subgroup $H \leq G$, the following diagram commutes:
$$\xymatrix{
N_{hG}^{Ad_G} \ar[r]^-{\tau_H^G} \ar[d]_\pi & N_{hH}^{Ad_H} \ar[d]^\pi \\
M_{hG}^{Ad_G + \nu} \ar[r]_-{\tau_H^G} & M_{hH}^{Ad_H + \nu}
}$$
Here $\pi$ are equivariant Pontrjagin-Thom collapse maps, as in Lemma \ref{ptmap_lemma}.  Since $\nu$ is equivariant, the transfer also commutes with the Thom isomorphism for $\nu$.

\end{lemma}

Let $\sO$ be the homology operad: $\sO(k) = \Sigma h_*(Cact_{hS^1})$, and let $A = \Sigma^{1-d} h_*(LM_{hS^1})$.  Showing that $\unrho$ makes $A$ into an $\sO$-algebra requires demonstrating that this diagram commutes:
$$\xymatrix{
\sO(k) \otimes \sO(l) \otimes A^{\otimes k+l-1} \ar[r]^-{1 \otimes \unrho} \ar[d]_{\uncirc_i \otimes 1} & \sO(k) \otimes A^{\otimes k} \ar[d]^{\unrho} \\
\sO(k+l-1) \otimes A^{\otimes k+l-1} \ar[r]_-{\unrho} & A
}$$
where $\unrho$, in the upper horizontal arrow, is the action on the $i^{\rm th}$ through $i+l^{\rm th}$ terms in $A^{\otimes k+l-1}$.

In proving the nonequivariant version, Cohen and Voronov's account \cite{cv} considers a space
$$C(k) \circ_i C(l)M \subseteq C(k) \times C(l) \times LM^{k+l-1}$$
This is the subspace of pairs of cacti $(c, d)$ and loops $(f_1, \dots, f_{k+l-1})$ which agree at the points of intersection of the composite cactus $c \circ_i d$, thus defining a map from it into $M$.  Note that this is a subspace of finite codimension.

The commutativity of the following diagram was established in \cite{cv}, Theorem 2.3.1.
$$\xymatrixcolsep{3pc}
\xymatrix{
C(k) \times C(l) \times LM^{k+l-1} \ar[dd]_{\circ_i \times 1} & C(k) \times C(l)M \times LM^{k-1}\ar[l]_-{\supseteq} \ar[r]^-{1 \times concat \times 1} & C(k) \times LM^{k} \\
 & C(k) \circ_i C(l)M \ar[dl]^{\circ_i \times 1} \ar[d]_{\circ_i \times 1} \ar[dr]_{concat} \ar[r]^-{partial}_-{concat} \ar[u]^{\subseteq} \ar[ul]^{\supseteq} & C(k)M \ar[u]_{\subseteq} \ar[d]_{concat} \\
C(k+l-1) \times LM^{k+l-1} & C(k+l-1)M \ar[l]^-{\supseteq} \ar[r]_-{concat} & LM
}$$
Here (partial) concatenations of loops in the upper right square of the diagram concatenate the $i^{\rm th}$ through $i+l^{\rm th}$ loops into the $i^{\rm th}$ factor of $LM^{\times k}$.  After replacing subspace inclusions with their associated umkehr maps, the iterated action ($\rho \circ \rho$) for $h_*(Cact)$ on $h_*(LM)$ is given along the top and right, and the action composed with substitution ($\rho \circ (\circ_i \times 1)$) is given along the bottom and left.  Consequently both are given by the map along the diagonal, and thus are equal.

To show associativity of the operad action in the equivariant case, we will need an action of $S^1$ on $Cact(k) \times Cact(l) \times LM^{k+l-1}$.  For $1 \leq i \leq k$, let the \emph{$i^{th}$ iterated homotopy diagonal action} be given by
$$\theta \cdot (c, d, f_1, \dots, f_{k+l-1}) = (\theta \cdot c, \Delta_i(c)(\theta) \cdot d, \Delta(c \circ_id)(\theta) \cdot (f_1, \dots, f_{k+l-1}))$$
Using similar arguments to the proof of Lemma \ref{action_lemma}, one can see that this is a group action.  It is evident that the vertical map $\circ_i \times 1$ is equivariant.  

We will show in Proposition \ref{coEnd_prop} that one can rewrite the $i^{\rm th}$ iterated homotopy diagonal action as
$$\theta \cdot (c, d, f_1, \dots, f_{k+l-1}) = (\theta \cdot c, \Delta_i(c)(\theta) \cdot d, (\Delta(c) \circ_i \Delta(d))(\theta) \cdot(f_1, \dots, f_{k+l-1}))$$
where $\Delta(c) \circ_i \Delta(d):S^1 \to (S^1)^{\times k+l-1}$ is the $i^{\rm th}$ operadic composition of the map $\Delta(c): S^1 \to (S^1)^{\times k}$ with the map $\Delta(d): S^1 \to (S^1)^{\times l}$.
%

Using this description of the action, we see that the subspaces in the upper middle and center positions of the diagram are equivariant subspaces, and that the various concatenation maps that originate from them in the diagram are equivariant.

Extend this to an action of $(S^1)^{\times k+l+1}$ on $C(k) \times C(l) \times LM^{\times k+l-1}$, where the first $S^1$ acts via the $i^{\rm th}$ iterated homotopy diagonal action, and the remaining $k+l$ factors act termwise on the last $k+l$ terms of $C(k) \times C(l) \times LM^{\times k+l-1}$.

\begin{lemma} \label{cactus_transfer_action_lemma}

The generalized string homology of $M$, $\Sigma^{1-d} h_*(LM_{hS^1})$, is an algebra over $\Sigma h_*(Cact_{hS^1})$ via the map $\unrho$.

\end{lemma}

\begin{proof}

Examine $\unrho \circ (\uncirc_i)$, the action composed with substitution in the equivariant case.  This is
$$\unrho \circ (\uncirc_i) = ((concat_{hS^1})_* \circ T \circ \pi_* \circ \tau_*) \circ (\circ_i \circ \tau_*)$$
(we suppress the equivalence $j$ in the definition of $\uncirc_i$).  The first transfer performed is with respect to the $i^{\rm th}$ homotopy diagonal action of $S^1 \leq S^1\times S^1$ on $C(k) \times C(l)$, and the second is with respect to the the homotopy diagonal action of $S^1 \leq (S^1)^{k+l}$ acting by the homotopy diagonal action on $C(k+l-1) \times LM^{\times k+l-1}$.  Since $\circ_i$ is $S^1$ equivariant, and transfers commute with equivariant maps, this can be rewritten as
$$(concat_{hS^1})_* \circ T \circ \pi_* \circ (\circ_i) \circ \tau_* \circ \tau_*$$
The iterated transfer $\tau_* \circ \tau_*$ may be identified with the transfer $\tau_*^{big}$ taken with respect to the subgroup $S^1 \leq (S^1)^{\times k+l+1}$ acting by the $i^{\rm th}$ iterated homotopy diagonal action.  So $\unrho \circ \uncirc_i$ is the composite $(concat_{hS^1})_* \circ T \circ \pi_* \circ (\circ_i) \circ \tau_*^{big}$.

Using the commutativity of the bottom left triangle in the diagram above, this is then
$$(concat_{hS^1})_* \circ T^{big} \circ \pi_*^{big} \circ \tau^{big}$$
Here $concat$, $T^{big}$, $\pi^{big}$ are the concatenation, Thom isomorphism, and Pontrjagin-Thom collapse with respect to the diagonal of the diagram:
$$
\xymatrix{(C(k) \times C(l) \times LM^{k+l-1})_{hS^1} & (C(k) \circ_i C(l)M)_{hS^1} \ar[r]^-{concat_{hS^1}} \ar[l]_-{\supseteq} & LM_{hS^1} }
$$

The iterated action is 
$$
\begin{array}{rcl}
\unrho \circ \unrho & = & ((concat_{hS^1})_* \circ T \circ \pi_* \circ \tau_*) \circ ((concat_{hS^1})_* \circ T \circ \pi_* \circ \tau_*) \\ 
 & = & (concat_{hS^1})_* \circ T \circ (\pi \circ concat_{hS^1})_* \circ \tau_* \circ T \circ \pi_* \circ \tau_* \\
 & = & (concat_{hS^1})_* \circ T \circ (\pi \circ concat_{hS^1})_* \circ T \circ \pi_* \circ \tau_* \circ \tau_*
\end{array} 
$$
The first step uses the fact that transfer maps commute with equivariant maps, and the second uses Lemma \ref{equivpt_lemma}.  As previously, we identify the composite of iterated transfer maps as the relative transfer $\tau^{big}$ with respect to the $i^{\rm th}$ iterated homotopy diagonal action; in this case we use the second description of the action given above.

Commutativity of the top right triangle of the commutative diagram allows us to conclude that
$$\unrho \circ \unrho  = (concat_{hS^1})_* \circ T^{big} \circ \pi_*^{big} \circ \tau^{big}$$
which gives the result. 

\end{proof}

The proof of the homotopy fixed point version (i.e., Corollary \ref{fixed_string_corollary}) is analogous, and left to the reader.

\subsection{The string bracket} \label{bracket_subsection}

The goal of this section is to relate the results of the previous section to the Chas-Sullivan string bracket and to construct a homotopy fixed point analogue. \\

\noindent {\it Proof of Theorem \ref{string_homology_theorem}.}  The bulk of the theorem follows from Lemma \ref{cactus_transfer_action_lemma} and Theorem \ref{same_at_last_thm}; we just need to identify the gravity algebra structure in terms of the operations $\overline{m}_k$.  By definition,
$$\overline{m}_k(a_1 \otimes \dots \otimes a_k) = p(\tau_*(a_1) \cdots \tau_*(a_k))$$
We would like to show that this is the same as the gravity operation induced by the generator $\iota_k := \{\cdot ,\ldots , \cdot \}$ of
$$H_1(Cact_{bS^1}(k)) = H_0(Cact(k)_{hS^1}) = \Z$$
(i.e., the point class $\iota_k \in H_0(Cact(k)_{hS^1}) \cong H_0(Cact(k))$).

Examine the definition of the gravity algebra structure on $\Sigma^{1-d} H_*(LM_{hS^1})$:
$$\unrho = (concat_{hS^1})_* \circ T \circ \pi_* \circ \tau_*$$
The map $\tau$ here is the relative transfer for $S^1 \times \{ 1 \}^{\times k} < (S^1)^{\times k+1}$; this is equal to the (global) transfer for the action of $\{ 1 \} \times (S^1)^{\times k}$.  So we compute
$$\{ a_1, \dots, a_k \} := \unrho(\iota_k ; a_1, \dots, a_k) = (concat_{hS^1})_* \circ T \circ \pi_*(\iota_k; \tau_*(a_1), \dots, \tau_*(a_k))$$

Nonequivariantly,
$$concat_* \circ T \circ \pi_*(\iota_k ; b_1, \dots, b_k) = b_1 \cdots b_k$$ 
That is, the point class of $H_0(Cact(k))$ induces the iterated loop product.  Thus  
$$(concat_{hS^1})_* \circ T \circ \pi_*(\iota_k; \tau_*(a_1), \dots, \tau_*(a_k)) = p(\tau_*(a_1) \cdots \tau_*(a_k))$$
\qed

Recall that for a naive $G$-spectrum $X$ one may define a \emph{cotransfer map}
$$c^G: S^{Ad_G} \sma X \to X^{hG}$$
(with properties akin to the transfer) as the composite
$$\xymatrix{
S^{Ad_G} \sma X \ar[r] & S^{Ad_G} \sma_{hG} X \ar[r]^-{n^G} & X^{hG}
}$$

We can define an analogue of the Chas-Sullivan bracket for the homotopy fixed points $H_*(LM^{hS^1})$ using the circle cotransfer $c = c^{S^1}$.

\begin{definition}

Define an operation 
$$\{ \cdot, \cdot \}: H_p(LM^{hS^1}) \otimes H_q(LM^{hS^1}) \to H_{p+q+1-d}(LM^{hS^1})$$
by $\{ a, b \} := c_*(i_*(a) \cdot i_*(b))$, where $\cdot$ is the loop product, and
$$i: LM^{hS^1} := (\Sigma^{\infty} LM_+)^{hS^1} \to \Sigma^{\infty} LM_+$$
is the (forgetful) inclusion of the homotopy fixed points.  We call $\{ \cdot, \cdot \}$ the \emph{homotopy fixed point string bracket}.  Using the same construction, one may define a bracket on $H^{c}_*(LM^{hS^1})$.

\end{definition}

As in the homotopy orbit case, we see that the action of the gravity operad induces this bracket:

\begin{lemma}

The homotopy fixed point string bracket is governed by the gravity algebra structure on $\Sigma^{1-d} H_*(LM^{hS^1})$; that is, it is induced by the generator of $H_1(\sC_2(2)^{hS^1}) \cong \Z$.  The same is true in continuous homology.

\end{lemma}

Apply Corollary \ref{splitting_corollary} to give the following result.

\begin{corollary} \label{fixed_lie_corollary}

The homotopy fixed point string bracket gives $A := H_*^{c}(LM^{hS^1})$ the structure of a Lie algebra.  It is also a module over the ring $H^{-*}(BS^1) = \Z[c]$, and the subspace $cA$ is an abelian Lie ideal.

\end{corollary}

The functoriality discussion in section \ref{homotopy_section} and Lemma \ref{op_map_lemma} can be adapted to the cactus operad to give a diagram of operads (up to homotopy):
$$\xymatrix{
Cact_{bS^1} \ar[rr]^-{n^{S^1}} \ar[dr]_-{\tau^{S^1}} & & Cact^{hS^1} \ar[dl]^-i \\
 & \Sigma^{\infty} Cact_+ &
}$$
Note that the two diagonal maps induce, upon application of $H_*$, the embedding of $Grav$ in $e_2$.  This makes $\Sigma^{-d} H_*(LM)$ into a gravity algebra; for instance the bracket $\{ \cdot, \cdot \}$ is given by the loop bracket $[ \cdot, \cdot ]$ (or Browder operation).  There is similarly a diagram of algebras:
$$\xymatrix{
LM_{bS^1} \ar[rr]^-{n^{S^1}} \ar[dr]_-{\tau^{S^1}} & & LM^{hS^1} \ar[dl]^-i \\
 & \Sigma^{\infty} LM_+ &
}$$
Thus, in homology, $n^{S^1}_*$, $\tau^{S^1}_*$, and $i_*$ are maps of gravity algebras.  One can see that $\tau^{S^1}_*$ is a map of Lie algebras directly; recall that for $\alpha, \beta \in H_*(LM)$,
$$[\alpha, \beta] = \Delta(\alpha \cdot \beta) - \Delta(\alpha) \cdot \beta - (-1)^{|\alpha|} \alpha \cdot \Delta(\beta)$$
So for $a, b \in H_*(LM_{hS^1})$,
$$\begin{array}{rcl}
[\tau_*^{S^1}(a), \tau_*^{S^1}(b)] & = & \Delta(\tau_*^{S^1}(a) \cdot \tau_*^{S^1}(b)) - \Delta(\tau_*^{S^1}(a)) \cdot \tau_*^{S^1}(b) - \\
 & & (-1)^{|a|+1} \tau_*^{S^1}(a) \cdot \Delta(\tau_*^{S^1}(b)) \\
 & = & \Delta(\tau_*^{S^1}(a) \cdot \tau_*^{S^1}(b))
\end{array}$$
The two latter terms vanish because $\Delta = \tau_*^{S^1} \circ \pi_*$ and $\pi_* \circ \tau_*^{S^1} = 0$.  The bracket in $H_*(LM_{hS^1})$ is defined by $\{ a, b \} = \pi_*(\tau_*^{S^1}(a) \cdot \tau_*^{S^1}(b))$, so 
$$\tau_*^{S^1} \{ a, b \} = \tau_*^{S^1} \pi_*(\tau_*^{S^1}(a) \cdot \tau_*^{S^1}(b)) = \Delta(\tau_*^{S^1}(a) \cdot \tau_*^{S^1}(b))$$
So $\tau_*^{S^1}$ is a map of Lie algebras.

\section{Proofs of technical results} \label{tech_section}

We collect the proofs of some of the more technical results in this section.

\subsection{The $S^1$ action on the configuration space of the plane}

It is well known (see, e.g., \cite{may}) that, if $F(X, k)$ denotes the ordered configuration space of $k$ points in a topological space $X$, the map
$$\sC_n(k) \to F(\R^n, k)$$
which carries a little disk to its center is an $SO(n)$-equivariant equivalence.  Our goal in this section is a proof of Lemma \ref{free_lemma}; the previous says that to prove this, we may as well replace $\sC_2(k)$ with $F(\C, k) = F(\R^2, k)$.

There is a fibration
$$p: F(\C, k) \to F(\C, 2)$$
given by $(x_1, \dots, x_k) \mapsto (x_1, x_2)$.  The fibre of $p$ over the point $(0, 1) \in F(\C, 2)$ is $F(\C \setminus \{ 0, 1 \}, k-2)$.  Denote by $i$ the inclusion of the fibre:
$$i: F(\C \setminus \{ 0, 1 \}, k-2) \subseteq F(\C, 2)$$
Let $\pi: F(\C, k) \to F(\C, k)/S^1$ be the quotient map.  Note that when $k \geq 2$, Lemma \ref{free_space_lemma} implies that $F(\C, k)/S^1 \simeq F(\C, k)_{hS^1}$.

\begin{proposition} \label{fibre_prop}

For $k \geq 2$, the composite $\pi \circ i: F(\C \setminus \{ 0, 1 \}, k-2) \to F(\C, k)/S^1$ is an equivalence.

\end{proposition}

\begin{proof}

The source and target of the fibration $p$ are both $S^1$-spaces, and $p$ is clearly equivariant.  Therefore $p/S^1$ is a fibration, and we get a map of fibrations
$$\xymatrix{
F(\C \setminus \{ 0, 1 \}, k-2) \ar[r]^-{i} \ar[d]_{=} & F(\C, k) \ar[r]^-{p} \ar[d]_{\pi} & F(\C, 2) \ar[d]_{\pi} \\
F(\C \setminus \{ 0, 1 \}, k-2) \ar[r]^-{i/S^1} & F(\C, k)/S^1 \ar[r]^-{p/S^1} & F(\C, 2)/S^1
}$$

However, $F(\C, 2)$ is $S^1$-equivariantly equivalent to $S^1$, so the base of the lower fibration is contractible.  The result follows from the long exact sequence in homotopy groups. 

\end{proof}

Write 
$$B = i_*(H_*(F(\C \setminus \{ 0, 1 \}, k-2))) \subseteq H_*(F(\C, k))$$
In proving Lemma \ref{free_lemma}, we will show that $B$ is a basis for the action of $\Lambda[\Delta]$ on $H_*(F(\C, 2))$.

\begin{corollary}

$\Delta(B) \cap B = 0$.

\end{corollary}

\begin{proof}

In general, $\pi_*(\Delta(x)) = 0$.  By the previous proposition, $\pi_* \circ i_*$ is an isomorphism, so the image of $i_*$ must intersect the image of $\Delta$ only in $0$. 

\end{proof}

\noindent {\it Proof of Lemma \ref{free_lemma}.}  To show that $H_*(F(\C, k))$ is free over $\Lambda[\Delta]$ on $B$, we must show that 
$$H_*(F(\C, k)) = B \oplus \Delta(B)$$
By the previous corollary, it suffices to show that $H_*(F(\C, k)) = B + \Delta(B)$.

If we examine the Serre spectral sequence for the fibration $p$, we obtain an $E_2$ term
$$E_2 = H_*(F(\C, 2)) \otimes H_*(F(\C \setminus \{ 0, 1 \}, k-2))$$
since the system of local coefficients is simple.  This spectral sequence must necessarily collapse, as any differential is of horizontal length at least $2$, and $E_2^{p, q} = 0$ for $p<0$ or $p>1$.  So $E_2$ is the associated graded group for a filtration of $H_*(F(\C, k))$, and $B$ corresponds to the subgroup 
$$B' = 1 \otimes H_*(F(\C \setminus \{ 0, 1 \}, k-2))$$

The map $p$ is an $S^1$-equivariant fibration, so for $b \in B$, represented by $1 \otimes b' \in B'$, up to terms in a higher filtration,
$$\Delta(b) = \Delta \otimes b'$$
Since elements of the form $1 \otimes b'$ and $\Delta \otimes b'$ generate $H_*(F(\C, k))$, the lemma follows.  \\ \qed

\begin{remark}

Notice that Proposition \ref{fibre_prop} gives an alternate computation of $e_2^{S^1}(k) = {e_2}_{S^1}(k)$ for $k \geq 2$:
$${e_2}_{S^1}(k) \cong \Sigma H_*(F(\C \setminus \{ 0, 1 \}, k-2))$$

This gives us an alternative proof of Proposition \ref{gravity_prop}, because there is an equivalence $\sM_{0, k+1} \cong F(\C \setminus \{ 0, 1 \}, k-2)$:  For a configuration of points $(x_1, \dots, x_{k+1})$ in $\C P^1$, there is a unique automorphism $\phi$ of $\C P^1$ carrying $(x_1, x_2, x_3)$ to $(0, 1, \infty)$.  So the map $\sM_{0, k+1} \to F(\C \setminus \{ 0, 1 \}, k-2)$ given by
$$(x_1, \dots, x_{k+1}) \mapsto (\phi(x_4), \dots, \phi(x_{k+1}))$$
is an equivalence, with obvious inverse.

\end{remark}

\subsection{Results for cacti}

The essential part of the proof of Lemma \ref{action_lemma} is the following:

\begin{lemma} \label{formula_lemma}

For $c \in Cact(k)$, $\theta, \phi \in S^1$ (whose product we will write multiplicatively), 
$$\Delta(c) (\phi \cdot \theta) = (\Delta(\theta\cdot c)(\phi))\cdot(\Delta(c)(\theta))$$
where $\cdot$ on the right side is termwise multiplication in $(S^1)^{\times k}$.

\end{lemma}

\begin{proof}

We claim that it is enough to show this in the limited case for all $\theta$, $c$, and when $\phi$ is small enough that one of the following two requirements holds:

\begin{enumerate}

\item $\phi$ does not move the outer marked point of $\theta \cdot c$ off the circle on which it already lies. \label{samecirc}

\item $\phi$ moves the outer marked point of $\theta \cdot c$ to the point of intersection of the circle that it lies on with the next circle of the cactus. \label{nextcirc}

\end{enumerate}

Every $\phi$ may be written $\phi = \phi_1 \cdot \phi_2 \cdots \phi_n$ for some choice of $n$, with $\phi_1$ of type \ref{samecirc}, and $\phi_2, \dots, \phi_n$ of type \ref{nextcirc}.  The lemma for general $\phi$ then follows by induction on the length $n$.

In general, suppose that $\Delta (c) (\theta) = (\theta_1, \dots, \theta_k)$ and that the marked point of $\theta \cdot c$ is in the $i^{\rm th}$ circle.  Let $\phi$ be of type \ref{samecirc}; then
$$\Delta (c) (\phi \cdot \theta) = (\theta_1, \dots, \theta_{i-1}, \phi \cdot \theta_i, \theta_{i+1}, \dots, \theta_k)$$
but by assumption, $\Delta(\theta \cdot c)(\phi) = (1, \dots, 1, \phi, 1, \dots, 1)$ (with $\phi$ in the $i^{\rm th}$ term).  So the lemma follows in this case.  If $\phi$ is of type \ref{nextcirc}, the same argument works (as it must, by continuity). 

\end{proof}

\noindent {\it Proof of Lemma \ref{action_lemma}.}  To show that the homotopy diagonal action is a group action, we must show that for $\theta, \phi \in S^1$, $c \in Cact(k)$, and $c_i \in Cact(n_i)$,
$$\phi \cdot (\theta \cdot (c; \; c_1, \dots, c_k)) = (\phi \cdot \theta) \cdot (c; \; c_1, \dots, c_k)$$
The right hand side can be written as
$$((\phi \cdot \theta) \cdot c; \; \Delta(c)(\phi \cdot \theta) \cdot (c_1, \dots, c_k))$$
and the left hand side is
$$(\phi \cdot (\theta \cdot c); \; (\Delta(\theta\cdot c)(\phi))\cdot(\Delta(c)(\theta)) \cdot (c_1, \dots, c_k))$$
so the first part of Lemma \ref{action_lemma} follows from Lemma \ref{formula_lemma}.

The homotopy diagonal action is shown to be homotopic to the diagonal action in the proof of Theorem 5.3.4 in \cite{kaufcacti}.  One can see that the homotopy $h_t$ that Kaufmann exhibits induces a group action for all $t$ by direct computation. \\ \qed

\noindent {\it Proof of Proposition \ref{equiv_prop}.}  It will be more transparent (and equivalent) to show that the full substitution map
$$\gamma: C(k) \times C(n_1) \times \dots \times C(n_k) \to C(\sum n_i)$$
is equivariant.  To see this, we need to verify two facts.  First, we must show that $\gamma(\theta \cdot(c_0; c_1, \dots, c_k))$ has the same ``shape'' as $\gamma(c_0; c_1, \dots, c_k)$ -- the cacti $c_1, \dots, c_k$ are attached to each other at the same points in the two substitutions.  Second, we show that the outer marked point on $\gamma(\theta \cdot(c_0; c_1, \dots, c_k))$ is that of $\gamma(c_0; c_1, \dots, c_k)$, advanced by $\theta$.

The substitution map glues $c_1, \dots, c_k$ together by identifying the inner marked point of the $i^{\rm th}$ circle of $c_0$ with the outer marked point of $c_i$ via the pinching $\nabla: S^1 \to c_i$.  The homotopy diagonal action advances the outer markings on the $c_i$ ($i>0$) in a piecewise fashion, by the same amount the outer marking on the $i^{\rm th}$ circle of $c_0$ is advanced.  Consequently the attaching points of the new cactus are identical.

That the outer marked point on $\gamma(\theta \cdot(c_0; c_1, \dots, c_k))$ is the same as the outer marked point on $\theta \cdot \gamma(c_0; c_1, \dots, c_k)$ is almost precisely the definition of the homotopy diagonal action. \\ \qed \\

Recall that, for any space $X$, one may define the \emph{co-endomorphism operad} $coEnd_X$ by
$$coEnd_X(k) := Map (X, X^{\times k})$$
whose substitution map is composition: for $f \in coEnd_X(k)$ with $f(\theta) = (\theta_1, \dots, \theta_k)$, define $\circ_i: coEnd_X(k) \times coEnd_X(l) \to coEnd_X(k+l-1)$ by
$$f \circ_i g (\theta) = (\theta_1, \dots, \theta_{i-1}, g(\theta_i), \theta_{i+1}, \dots, \theta_k)$$
Equivalently, $\gamma(f, g_1, \dots, g_k) = (g_1 \times \dots \times g_k) \circ f$.

The action of $\Sigma_k$ is defined by postcomposition with the action on the target $X^{\times k}$.

\begin{proposition} \label{coEnd_prop}

The homotopy diagonals $\Delta: Cact(k) \to Map(S^1, (S^1)^{\times k})$ define a map of operads
$$Cact \to coEnd_{S^1}$$

\end{proposition}

\begin{proof}

Equivalently, we must show that if $\Delta(c)(\theta) = (\theta_1, \dots, \theta_k)$ for $c \in Cact(k)$, then
$$\Delta(c \circ_i d)(\theta) = (\theta_1, \dots, \theta_{i-1}, \Delta(d)(\theta_i), \theta_{i+1}, \dots, \theta_k)$$
This is precisely how $\Delta$ is defined. 

\end{proof}


\thanks{\noindent {\bf Acknowledgements.}  It is a pleasure to thank the host of people whose helpful comments have aided the evolution of this project: Andrew Blumberg, Michael Ching, Ralph Cohen, Tom Fiore, John Greenlees, John Klein, Johann Leida, Jack Morava, Peter May, Holger Reich, and John Rognes.  I would particularly like to thank Ezra Getzler for alerting me to the connection between the gravity operad and transfer operad for the circle action on the little disk operad.  The referee also deserves thanks for several suggestions which improved the exposition of the paper.  This paper was written at the University of Wisconsin and MSRI; I would like to thank both institutions for their support.}


\bibliography{biblio}

\providecommand{\bysame}{\leavevmode\hbox to3em{\hrulefill}\thinspace}
\providecommand{\MR}{\relax\ifhmode\unskip\space\fi MR }
\providecommand{\MRhref}[2]{%
  \href{http://www.ams.org/mathscinet-getitem?mr=#1}{#2}
}
\providecommand{\href}[2]{#2}
\begin{thebibliography}{EKMM97}

\bibitem[ACD89]{acd}
A.~Adem, R.~L. Cohen, and W.~G. Dwyer, \emph{Generalized {T}ate homology,
  homotopy fixed points and the transfer}, Algebraic topology (Evanston, IL,
  1988), Contemp. Math., vol.~96, Amer. Math. Soc., Providence, RI, 1989,
  pp.~1--13.

\bibitem[ACG05]{acg}
H.~Abbaspour, R.~L. Cohen, and K.~Gruher, \emph{{String topology of Poincare
  duality groups}}, preprint: math.AT/0511181 (2005).

\bibitem[BR05]{bruner_rognes}
R.~R. Bruner and J.~Rognes, \emph{Differentials in the homological homotopy
  fixed point spectral sequence}, Algebr. Geom. Topol. \textbf{5} (2005),
  653--690 (electronic).

\bibitem[BV73]{bv}
J.~M. Boardman and R.~M. Vogt, \emph{Homotopy invariant algebraic structures on
  topological spaces}, Lecture Notes in Mathematics 347, Springer-Verlag,
  Berlin, 1973.

\bibitem[Car91]{carlssonloop}
G.~Carlsson, \emph{On the homotopy fixed point problem for free loop spaces and
  other function complexes}, $K$-Theory \textbf{4} (1991), no.~4, 339--361.

\bibitem[Chi05]{ching}
Michael Ching, \emph{Bar constructions for topological operads and the
  {G}oodwillie derivatives of the identity}, Geom. Topol. \textbf{9} (2005),
  833--933 (electronic).

\bibitem[CHV06]{cv}
Ralph~L. Cohen, Kathryn Hess, and Alexander~A. Voronov, \emph{String topology
  and cyclic homology}, Advanced Courses in Mathematics. CRM Barcelona,
  Birkh\"auser Verlag, Basel, 2006, Lectures from the Summer School held in
  Almer\'\i a, September 16--20, 2003.

\bibitem[CJ02]{cj}
R.~L. Cohen and J.~D.~S. Jones, \emph{A homotopy theoretic realization of
  string topology}, Mathematische Annalen \textbf{324 no. 4} (2002), 773--798.

\bibitem[CLM76]{clm}
F.~R. Cohen, T.~J. Lada, and J.~P. May, \emph{The homology of iterated loop
  spaces}, Springer-Verlag, Berlin, 1976, Lecture Notes in Mathematics, Vol.
  533.

\bibitem[CS01]{cs}
M.~Chas and D.~Sullivan, \emph{String topology}, preprint: math.GT/9911159
  (2001).

\bibitem[Dev99]{devadoss}
Satyan~L. Devadoss, \emph{Tessellations of moduli spaces and the mosaic
  operad}, Homotopy invariant algebraic structures (Baltimore, MD, 1998),
  Contemp. Math., vol. 239, Amer. Math. Soc., Providence, RI, 1999,
  pp.~91--114.

\bibitem[EHKR05]{ehkr}
P.~Etingof, A.~Henriques, J.~Kamnitzer, and E.~Rains, \emph{The cohomology ring
  of the real locus of the moduli space of stable curves of genus 0 with marked
  points}, preprint: math.AT/0507514 (2005).

\bibitem[EKMM97]{EKMM}
A.~D. Elmendorf, I.~Kriz, M.~A. Mandell, and J.~P. May, \emph{Rings, modules,
  and algebras in stable homotopy theory}, Mathematical Surveys and Monographs,
  vol.~47, American Mathematical Society, Providence, RI, 1997, With an
  appendix by M. Cole.

\bibitem[EM97]{elmay}
A.~D. Elmendorf and J.~P. May, \emph{Algebras over equivariant sphere spectra},
  J. Pure Appl. Algebra \textbf{116} (1997), no.~1-3, 139--149, Special volume
  on the occasion of the 60th birthday of Professor Peter J. Freyd.

\bibitem[Fie98]{fied}
Z.~Fiedorowicz, \emph{Constructions of ${E}_n$ operads}, preprint:
  math.AT/9808089 (1998).

\bibitem[Get94]{getz2d}
E.~Getzler, \emph{Two-dimensional topological gravity and equivariant
  cohomology}, Comm. Math. Phys. \textbf{163} (1994), no.~3, 473--489.

\bibitem[Get95]{getz0}
\bysame, \emph{Operads and moduli spaces of genus {$0$} {R}iemann surfaces},
  The moduli space of curves (Texel Island, 1994), Progr. Math., vol. 129,
  Birkh\"auser Boston, Boston, MA, 1995, pp.~199--230.

\bibitem[GJ94]{gj}
E.~Getzler and J.D.S. Jones, \emph{Operads, homotopy algebra and iterated
  integrals for double loop spaces}, preprint: hep-th/9403055 (1994).

\bibitem[GK94]{gk}
V.~Ginzburg and M.~Kapranov, \emph{Koszul duality for operads}, Duke Math. J.
  \textbf{76} (1994), no.~1, 203--272.

\bibitem[GM95]{gm}
J.~P.~C. Greenlees and J.~P. May, \emph{Generalized {T}ate cohomology}, Mem.
  Amer. Math. Soc. \textbf{113} (1995), no.~543, viii+178.

\bibitem[GS06]{grusalv}
K.~Gruher and P.~Salvatore, \emph{{Generalized string topology operations}},
  preprint: math.AT/0602210 (2006).

\bibitem[GW07]{gw}
K.~Gruher and C.~Westerland, \emph{{String topology prospectra and Hochschild
  cohomology}}, in preparation (2007).

\bibitem[Kau05]{kaufcacti}
R.~M. Kaufmann, \emph{On several varieties of cacti and their relations},
  Algebr. Geom. Topol. \textbf{5} (2005), 237--300 (electronic).

\bibitem[Kel05]{kelly}
G.~M. Kelly, \emph{On the operads of {J}. {P}. {M}ay}, Repr. Theory Appl.
  Categ. (2005), no.~13, 1--13 (electronic).

\bibitem[Kle01]{klein_dual}
J.~R. Klein, \emph{The dualizing spectrum of a topological group}, Math. Ann.
  \textbf{319} (2001), no.~3, 421--456.

\bibitem[KM94]{km}
M.~Kontsevich and Yu. Manin, \emph{Gromov-{W}itten classes, quantum cohomology,
  and enumerative geometry}, Comm. Math. Phys. \textbf{164} (1994), no.~3,
  525--562.

\bibitem[LMSM86]{lms}
L.~G. Lewis, Jr., J.~P. May, M.~Steinberger, and J.~E. McClure,
  \emph{Equivariant stable homotopy theory}, Lecture Notes in Mathematics, vol.
  1213, Springer-Verlag, Berlin, 1986, With contributions by J. E. McClure.

\bibitem[May72]{may}
J.~P. May, \emph{The geometry of iterated loop spaces}, Lecture Notes in
  Mathematics 271, Springer Verlag, Berlin, 1972.

\bibitem[May97]{may_defn}
J.~P. May, \emph{Definitions: operads, algebras and modules}, Operads:
  Proceedings of Renaissance Conferences (Hartford, CT/Luminy, 1995)
  (Providence, RI), Contemp. Math., vol. 202, Amer. Math. Soc., 1997, pp.~1--7.

\bibitem[MM02]{manmay}
M.~A. Mandell and J.~P. May, \emph{Equivariant orthogonal spectra and
  {$S$}-modules}, Mem. Amer. Math. Soc. \textbf{159} (2002), no.~755, x+108.

\bibitem[MS00]{madsch}
I.~Madsen and C.~Schlichtkrull, \emph{The circle transfer and {$K$}-theory},
  Geometry and topology: Aarhus (1998), Contemp. Math., vol. 258, Amer. Math.
  Soc., Providence, RI, 2000, pp.~307--328.

\bibitem[Rog05]{rognes}
J.~Rognes, \emph{Stably dualizable groups}, preprint: math.AT/0502184 (2005).

\bibitem[Str00]{strick}
N.~P. Strickland, \emph{{$K(N)$}-local duality for finite groups and
  groupoids}, Topology \textbf{39} (2000), no.~4, 733--772.

\bibitem[SW03]{salvwahl}
P.~Salvatore and N.~Wahl, \emph{Framed discs operads and {B}atalin-{V}ilkovisky
  algebras}, Q. J. Math. \textbf{54} (2003), no.~2, 213--231.

\bibitem[Vor05]{voruniv}
A.~A. Voronov, \emph{Notes on universal algebra}, Graphs and patterns in
  mathematics and theoretical physics, Proc. Sympos. Pure Math., vol.~73, Amer.
  Math. Soc., Providence, RI, 2005, pp.~81--103.

\bibitem[Wes07]{wes_string}
C.~Westerland, \emph{{String homology of spheres and projective spaces}},
  Algebr. Geom. Topol. \textbf{7} (2007), 309--325 (electronic).

\end{thebibliography}

\end{document}